\documentclass[graybox]{svmult}

\usepackage{type1cm}        
\usepackage{makeidx}         
\usepackage{graphicx}        
\usepackage{multicol}        
\usepackage[bottom]{footmisc}

\usepackage{newtxtext}       %
\usepackage[varvw]{newtxmath}       

\makeindex             

\usepackage{amsmath}
\usepackage{bigints}
\usepackage{mathtools}
\usepackage{rotating}
\usepackage{comment}

\newcommand{\bx}{\mathbf{x}}
\newcommand{\bX}{\mathbf{X}}
\newcommand{\by}{\mathbf{y}}
\newcommand{\EE}{\mathbb{E}}      

\newcommand{\R}{\mathbb{R}}     
\newcommand{\HH}{\mathcal{H}}
\newcommand{\X}{\mathcal{X}}
\renewcommand{\S}{\mathcal{S}}

\DeclareMathOperator*{\Argmax}{Arg\,max}
\DeclareMathOperator*{\Argmin}{Arg\,min}
\newcommand{\pyvar}[1]{\texttt{#1}}

\newcommand{\Kbarbar}{\overline{K}_{|m}}
\newcommand{\ISE}{\mathrm{ISE}}

\def\dd{\mbox{\rm d}}
\def\TT{^\top}
\def\0b{\mathbf{0}}
\def\Db{\mathbf{D}}
\def\kb{\mathbf{k}}
\def\Kb{\mathbf{K}}
\def\Ib{\mathbf{I}}
\def\Sb{\mathbf{S}}
\def\GP{\mathsf{GP}}
\newcommand{\bea}{\begin{eqnarray*}}
\newcommand{\eea}{\end{eqnarray*}}
\newcommand{\be}{\begin{eqnarray}}
\newcommand{\ee}{\end{eqnarray}}
\def\wb{\mathbf{w}}
\def\SE{{\mathscr E}}
\def\SN{{\mathscr N}}
\def\pb{\mathbf{p}}
\newcommand{\Kbarbarb}{\overline{\Kb}_{|m}}
\def\mve{\varepsilon}
\newcommand{\Kbarbarp}{\overline{K'}_{|m}}
\newcommand{\Kbarbarpb}{\overline{\Kb'}_{|m}}
\newcommand{\fin} {\mbox{}~\hfill{\lower-0.3ex\hbox{$\triangleleft$}}}
\def\e1{\mathrm{e}}
\newcommand{\bZ}{\mathbf{Z}}
\newcommand{\etab}{\mathbf{\eta}}

\begin{document}

\title*{Model predictivity assessment: incremental test-set selection and accuracy evaluation}
\titlerunning{Incremental test-set selection for model predictivity assessment}

\author{Elias Fekhari \and Bertrand Iooss \and Joseph Mur\'e \and Luc Pronzato \and Maria-Jo\~ao Rendas}
\authorrunning{E. Fekhari, B. Iooss, J. Mur\'e, L. Pronzato, M-J. Rendas}
\institute{
Elias Fekhari \and Bertrand Iooss$^\star$ \and Joseph Mur\'e \at EDF R\&D, 6 Quai Watier, 78401 Chatou, France - 
\email{elias.fekhari@edf.fr, bertrand.iooss@edf.fr, joseph.mure@edf.fr} - $^\star$ Corresponding author, phone: +33130877969
\and Luc Pronzato \and Maria-Jo\~ao Rendas \at CNRS, Universit\'e C\^ote d'Azur, Laboratoire I3S, B\^at. Euclide, Les Algorithmes, 2000 route des Lucioles, 06900 Sophia Antipolis cedex, France - \email{luc.pronzato@i3s.unice.fr, rendas@i3s.unice.fr}
}
%
%

\maketitle

\abstract*{
Unbiased assessment of the predictivity of models learnt by supervised machine learning (ML) methods  requires knowledge of the learned function over a reserved test set (not used by the learning algorithm). The quality of the assessment depends, naturally, on the properties of the test set and on the error statistic used to estimate the prediction error. In this work we tackle both issues, proposing a new predictivity criterion that carefully weights the individual observed errors to obtain a global error estimate, and using incremental experimental design methods to ``optimally''  select the test points on which the criterion is computed. Several incremental constructions are studied, including greedy-packing (coffee-house design), support points and kernel herding techniques. Our results show that the incremental and weighted versions of the  latter two, based on Maximum Mean Discrepancy concepts, yield superior performance.
An industrial test case provided by the historical French electricity supplier (EDF) illustrates the practical relevance of the methodology, indicating that it is an efficient alternative to expensive cross-validation techniques.
}
\abstract{
Unbiased assessment of the predictivity of models learnt by supervised machine learning (ML) methods  requires knowledge of the learned function over a reserved test set (not used by the learning algorithm). The quality of the assessment depends, naturally, on the properties of the test set and on the error statistic used to estimate the prediction error. In this work we tackle both issues, proposing a new predictivity criterion that carefully weights the individual observed errors to obtain a global error estimate, and using incremental experimental design methods to ``optimally''  select the test points on which the criterion is computed. Several incremental constructions are studied, including greedy-packing (coffee-house design), support points and kernel herding techniques. Our results show that the incremental and weighted versions of the  latter two, based on Maximum Mean Discrepancy concepts, yield superior performance.
An industrial test case provided by the historical French electricity supplier (EDF) illustrates the practical relevance of the methodology, indicating that it is an efficient alternative to expensive cross-validation techniques.
}
\vspace{.5cm}





\vspace{-0.3cm}
\section{Introduction}
\label{sec:1}

The development of tools for automatic diagnosis relying  on learned models 
imposes strict requirements on   model validation. 
For example, in industrial non-destructive testing (e.g. for the aeronautic or the nuclear industry), generalized automated inspection, which increases efficiency and lowers costs, must provide high performance guarantees \cite{eni19,hawpat21}. 
Establishing these guarantees requires availability of a 
reserved test set, i.e.\ a  data set
that has not been used either to train or to select the machine learning (ML) model \cite{borjir12,xugoo18,ioo21}.
Using the prediction residuals on this test set, an independent evaluation of the proposed ML model can be done, enabling the estimation of relevant performance metrics, such as the mean-squared error for regression problems, or the misclassification rate for classification problems.

The same need for independent test sets arises in the area of computer experiments, where computationally expensive simulation codes are often advantageously replaced 
by ML models, called surrogate models (or metamodels) in this context \cite{sanwil03,fanli06}.
Such surrogate models can be used, for instance, to  estimate the region of the input space that maps to specific values of the model outputs \cite{chebec14} with a significantly decreased computational load when compared to direct use of the original simulation code. 
Validation of these  surrogate models consists in estimating their predictivity, and can either rely on a suitably selected validation sample, or be done by cross-validation \cite{klesar00,ioobou10,demioo21}. 
One of the numerical studies presented in this paper will address an example of this situation of practical industrial interest in the domain of nuclear safety assessment, concerning the simulation of thermal-hydraulic phenomena inside nuclear pressurized water reactors, for which finely validated surrogate models have demonstrated their usefulness  \cite{lorzan11,marioo21}.

In this paper, we present methods to choose a ``good'' test set,
either within a given dataset or within the input space of the model,
as recently motivated in \cite{ioo21,josvak22}.
A first choice  concerns the size of the test set.
No optimal choice exists, and, when only a finite dataset is available,  classical ML handbooks \cite{hastib09,gooben16} provide different heuristics on how to split it, e.g., $80\% / 20\%$ between the training and test samples, or $50\% / 25\% / 25\%$ between the training, validation (used for model selection) and test samples.
We shall not formally address this point here (see \cite{xugoo18} for a numerical study of this issue),
but in the industrial case-study mentioned above we do study the impact of the ratio between the sizes of the training and test sets on the ability of assessing the quality of the surrogate model.
A second issue
concerns how the test sample is picked within the input space.
The simplest -- and most common  -- way to build a test sample is to extract an independent Monte Carlo sample \cite{hastib09}.
For small test sets, these randomly chosen points may fall too close to the training points or leave large areas of the input space unsampled, and a more constructive method to select points inside the input domain is therefore
preferable. 
Similar concerns motivate the use of space-filling designs when choosing a small set of runs for cpu-time expensive computer experiments on which a model will be identified \cite{fanli06,promul12}.



When the test set must be a subset of an initial dataset, the problem 
amounts to selecting a certain number of points within a finite collection of points. 
A review of classical methods for solving this issue is given in \cite{borjir12}.
For example, the CADEX and DUPLEX algorithms \cite{kensto69,sne77} 
can sequentially extract points from a database to include them in a test sample, using an inter-point distance criterion.
In ML, identifying within the dataset ``prototypes'' (set of data instances representative of the whole data set) and ``criticisms'' (data instances poorly represented by the prototypes) has recently been proposed to help model interpretation \cite{mol19}; 
the extraction of prototypes and criticisms relies on a Maximum Mean Discrepancy (MMD) criterion \cite{smogre07} (see e.g.\ \cite{prozhi20}, and \cite{pro21} for a performance analysis of greedy algorithms for MMD minimization). 

Several algorithms have also been proposed for the case where points need to be added to an already existing training sample.
When the goal is to  assess the quality of a model learnt using a known training set, one may be tempted to locate the test points the furthest away from the training samples, such that, in some sense, the union of the training and test sets is space-filling. As this paper shows, test sets built in this manner do enable a good assessment of the quality of models learnt with the training set if the observed residuals are appropriately weighted.
Moreover, the incremental augmentation of a design can be useful 
when the assessed model turns out to be of poor quality, or when an additional computational budget is available after a first study~\cite{sheraz17,shaapl21}.
Different empirical strategies have been proposed for incremental space-filling design \cite{ioobou10,crolae11,lilu17}, which basically entail the addition of new points in the zones poorly covered by the current design.
Shang and Apley~\cite{shaapl21} have recently proposed an improvement of the CADEX algorithm, called the Fully-sequential space-filling (FSSF) design; see also \cite{nogpro21} for an alternative version of coffee-house design enforcing boundary avoidance. 
Although they are developed for different purposes, nested space filling designs~\cite{qiaai09} and sliced space filling designs~\cite{qiawu09} can also be used to build sequential designs. 

In this work, we provide new insights into these subjects in two main directions: (\textit{i}) definition of new predictivity criteria through an optimal weighting of the test points residuals, and (\textit{ii}) use of test sets built by incremental space-filling algorithms, namely FSSF, support points \cite{makjos18} and kernel herding \cite{chewel10}, the latter two algorithms being typically used to provide a representative sample of a desired theoretical or empirical distribution.
Besides, this paper presents a numerical benchmark analysis comparing the behaviour of the three algorithms on a selected set of test cases.

This paper is organized as follows.
Section \ref{sec:3} defines the predictivity criterion considered and proposes different methods for its estimation.
Section \ref{sec:2} 
presents
the three algorithms used for test-point selection: 
FSSF, support points and kernel herding.
Our numerical results are presented in Sections \ref{sec:4} and \ref{sec:5}: in Section~\ref{sec:4} a test set is freely chosen within the entire input space, while in Section \ref{sec:5} an existing data set can be split into a training sample and a test set. Finally, Section \ref{sec:6} concludes 
and outlines some perspectives.

\section{Predictivity assessment criteria for an ML model}
\label{sec:3}

In this section, we propose a new criterion to assess the predictive performance of a model, derived from a standard model quality metric by suitably weighting the errors observed on the test set.
We denote by $\X \subset \R^d$ the space of the input variables $\bx=(x_1,\ldots,x_d)$ of the model.
Then let $y(\bx) \in \R$ (resp. $y(\bx') \in \R$) be the observed output at point $\bx \in \X$ (resp. $\bx' \in \X$). 
We denote by $(\bX_m,\by_m)$ the training sample, with $\by_m = [y(\bx^{(1)}),\ldots,y(\bx^{(m)})]^\top$.
The test sample is denoted by $(\bX_n,\by_n) = (\bx^{(m+i)},y(\bx^{(m+i)}))_{1\leq i \leq n}$.


\subsection{The predictivity coefficient}

Let $\eta_m(\bx)$ denote the prediction at point $\bx$ of a model 
 learned using $(\bX_m,\by_m)$ \cite{hastib09,raswil06}.
A classical measure for assessing the predictive ability of $\eta_m$, in order to evaluate its validity, is the predictivity coefficient.
Let $\mu$ denote the measure that weights how comparatively important it is to accurately predict $y$ over the different regions of $\X$. For example the input could be a random vector with known distribution: in that case, this distribution would be a reasonable choice for $\mu$.
The true (ideal) value of the predictivity is defined as the following normalization of the Integrated Square Error (ISE):
\begin{equation}\label{eq:Q2th}
Q_{\mathrm{ideal}}^2(\mu) = 1 - \frac{\ISE_\mu(\bX_m,\by_m)}{V_\mu} \,, 
\end{equation}
where
\begin{align*}
    \ISE_\mu(\bX_m,\by_m) &= \int_\X [y(\bx)-\eta_m(\bx)]^2 \, \dd\mu(\bx) \,, \\
    V_\mu &= \int_\X \left[y(\bx)-\int_\X y(\bx') d\mu(\bx') \right]^2 \, \dd\mu(\bx)\,.
\end{align*}
The ideal predictivity $Q_{\mathrm{ideal}}^2(\mu)$ is usually estimated by its empirical version calculated over the test sample $(\bX_n,\by_n)$, see \cite[p.~32]{davgam21}: 
\begin{equation}\label{eq:Q2test}
\widehat Q^2_n = 1 - \frac{ \sum_{\bx \in \bX_n}  \left[ y(\bx)-\eta_m(\bx)\right]^2}{\sum_{\bx \in \bX_n}  \left[y(\bx)-\overline{y}_n\right]^2}\,,
\end{equation}
where $\overline{y}_n=(1/n)\,\sum_{i=1}^n y(\bx^{(m+i)})$ denotes the empirical mean of the observations in the test sample. 
Note that the calculation of $\widehat Q^2_n$ only requires access to the predictor $\eta_m(\cdot)$. To compute $\widehat Q^2_n$, one does not need to know the training set which was used to build $\eta_m(\cdot)$.
$\widehat Q^2_n$ is the coefficient of determination (a standard notion in parametric regression) 
common in prediction studies \cite{klesar00,ioobou10}, 
often called ``Nash-Sutcliffe criterion'' \cite{NashS70}: 
it compares the prediction errors obtained with the model $\eta_m$ with those obtained when prediction equals the empirical mean of the observations. Thus, the closer $\widehat Q^2_n$ is to one, the more accurate the surrogate model is (for the test set considered). On the contrary, $\widehat Q^2_n$ close to zero (negative values are possible too) indicates poor predictions abilities, as there is little improvement compared to prediction by the simple empirical mean of the observations.  
The next section shows how a suitable weighting of the residual on the training sample may be key to improving the  estimation of $\widehat Q^2_n$. 

\subsection{Weighting the test sample}\label{S:weighting}



Let $\xi_n = n^{-1} \sum_{i=1}^n \delta_{\bx^{(i)}}$ be the empirical distribution of the prediction error, with $\delta_{\bx}$ the Dirac measure at $\bx$. In $\widehat Q^2_n$,  $\ISE_\mu(\bX_m,\by_m)$  is estimated by the empirical average of the squared residuals
$$
\ISE_{\xi_n}(\bX_m,\by_m) = \frac1n\, \sum_{i=1}^n  \left[ y(\bx^{(m+i)})-\eta_m(\bx^{(m+i)})\right]^2 \,.
$$
When the points $\bx^{(m+i)}$ of the test set $\bX_n$ are distant from the points of the training set $\bX_m$, the squared prediction errors $|y(\bx^{(m+i)})-\eta_m(\bx^{(m+i)})|^2$ tend to represent the worst possible error situations, and $\ISE_{\xi_n}(\bX_m,\by_m)$ tends to overestimate $\ISE_\mu(\bX_m,\by_m)$. In this section, we postulate a statistical model for the prediction errors
in order to be able to quantify this potential bias when sampling the residual process, enabling its subsequent correction .

In \cite{PR2021a}, the authors propose a weighting scheme for the test set when the ML model interpolates the train set observations. They suggest several variants corresponding to different constraints on the weights (e.g., non-negativity, summing to one). In the following, we consider the unconstrained version only, which in our experience works best. Let $\delta_m(\bx)=y(\bx)-\eta_m(\bx)$ denote the predictor error and assume it is  a realization of a Gaussian Process (GP) 
with zero mean and covariance kernel $\sigma^2\,K_{|m}$, 
which we shall note $\delta_m(\bx) \sim \GP(0,\sigma^2\,K_{|m})$, with
$$
\sigma^2\, K_{|m}(\bx,\bx')=\EE\{\delta_m(\bx)\delta_m(\bx')\}=\sigma^2\,\left[K(\bx,\bx')-\kb_m\TT(\bx)\Kb_m^{-1}\kb_m(\bx')\right]\,.
$$
Here, $\kb_m(\bx)$ denotes the column vector $[K(\bx,\bx^{(1)}),\ldots,K(\bx,\bx^{(m)})]\TT$ and $\Kb_m$ is the $m\times m$ matrix whose element $(i,j)$ is given by $\{\Kb_m\}_{i,j}=K(\bx^{(i)},\bx^{(j)})$, with $K$ a positive definite kernel. 
The rationale for using this model is simple: assume first a prior GP model $\GP(0,\sigma^2\,K)$ for the error process $\delta(\bx)$; if $\eta_m$ interpolates the observations $\by_m$, the errors observed at the design points $\bx^{(i)}$ equal zero, $i=1,\ldots,m$, leading finally to the posterior $\GP(0,\sigma^2\,K_{|m})$ for $\delta_m(\bx)$. 

However, the predictor $\eta_m$ is not always an interpolator, see Section~\ref{sec:5} for an example, so we extend the approach of \cite{PR2021a} to the general situation where $\eta_m$ does not necessarily interpolate the training data $\by_m$. 
The same prior $\GP(0,\sigma^2\,K)$ for $\delta(\bx)$ yields $\delta_m(\bx)\sim\GP(\widehat\delta_m(\bx),\sigma^2\,K_{|m})$, where 
\be\label{eq:deltahat}
\widehat\delta_m(\bx)=\kb_m\TT(\bx)\Kb_m^{-1}(\by_m-\etab_m)
\ee
is the Kriging interpolator for the errors, with $\etab_m=[\eta_m(\bx^{(1)}),\ldots,\eta_m(\bx^{(m)})]\TT$.

The model above allows us to study how well $\ISE_\mu(\bX_m,\by_m)$ is estimated using a given test set. Denote by $\overline{\Delta}^2(\xi_n,\mu;\bX_m,\by_m)$ the expected squared error when estimating $\ISE_\mu(\bX_m,\by_m)$ by $\ISE_{\xi_n}(\bX_m,\by_m)$,
\bea
\overline{\Delta}^2(\xi_n,\mu;\bX_m,\by_m) &=& \EE \left\{ \left[ \ISE_{\xi_n}(\bX_m,\by_m) - \ISE_\mu(\bX_m,\by_m) \right]^2 \right\}\\
&=& \EE\left\{ \left[\int_\X  \delta_m^2(\bx)\,\dd(\xi_n-\mu)(\bx)\right]^2\right\} \\
&=& \EE\left\{ \int_{\X^2}  \delta_m^2(\bx)\delta_m^2(\bx')\,\dd(\xi_n-\mu)(\bx)\dd(\xi_n-\mu)(\bx')\right\} \,.
\eea
Tonelli's theorem gives
\bea
\overline{\Delta}^2(\xi_n,\mu;\bX_m,\by_m)&=& \int_{\X^2}  \EE\{\delta_m^2(\bx)\delta_m^2(\bx')\}\,\dd(\xi_n-\mu)(\bx)\dd(\xi_n-\mu)(\bx') \,.
\eea
Since $\EE \left\{U^2 V^2 \right\} = 2\, \left(\EE\{UV\}\right)^2 + \EE \left\{ U^2 \right\} \EE \left\{ V^2 \right\}$
for any one-dimensional normal centered random variables $U$ and $V$, when $\eta_m(\bx)$ interpolates $\by_m$, we obtain
\be\label{Delta2}
\overline{\Delta}^2(\xi_n,\mu;\bX_m,\by_m)&=& \sigma^2\, d_{\Kbarbar}^2(\xi_n,\mu) \,,
\ee
where 
\be
    d_{\Kbarbar}^2(\xi_n,\mu) &=& \int_{\X^2}  \Kbarbar(\bx, \bx')\,\dd(\xi_n-\mu)(\bx)\dd(\xi_n-\mu)(\bx') \,, \nonumber\\
    \Kbarbar(\bx, \bx') &=& 2\, K_{|m}^2(\bx,\bx') + K_{|m}(\bx,\bx) K_{|m}(\bx',\bx')\,,  \label{Eq:Kbarbar}
\ee
and  we recognize $d_{\Kbarbar}^2(\xi_n,\mu)$ as the squared Maximum-Mean-Discrepancy (MMD)
between $\xi_n$ and $\mu$ for the kernel $\Kbarbar$; see \eqref{eq:MMD1} in Appendix~A. Note that $\sigma^2$ only appears as a multiplying factor in \eqref{Delta2}, with the consequence that $\sigma^2$ does not impact the choice of a suitable $\xi_n$.

When $\eta_m$ does not interpolate $\by_m$, 
similar developments still give
$\overline{\Delta}^2(\xi_n,\mu;\bX_m,\by_m) = \sigma^2\,d_{\Kbarbar}^2(\xi_n,\mu)$, with now
\bea
    \Kbarbar(\bx, \bx') &=& 2\, \left[K_{|m}(\bx,\bx') + 2\,\widehat\delta_m(\bx)\widehat\delta_m(\bx')\right]K_{|m}(\bx,\bx')  \\
    &&+\,\left[\widehat\delta_m^2(\bx)+K_{|m}(\bx,\bx)\right]\left[\widehat\delta_m^2(\bx')+K_{|m}(\bx',\bx')\right] \,,
\eea
where $\widehat\delta_m(\bx)$ is given by \eqref{eq:deltahat}.

The idea is to replace $\xi_n$, uniform on $\bX_n$, by a nonuniform measure $\zeta_n$ supported on $\bX_n$, $\zeta_n = \sum_{i=1}^n w_i \delta_{\bx^{(m+i)}}$ with weights $\wb_n = (w_1,\ldots,w_n)\TT$ chosen such that the estimation error $\overline{\Delta}^2(\xi_n,\mu;\bX_m,\by_m)$, and thus  $d_{\Kbarbar}^2(\zeta_n,\mu)$, is minimized. Direct calculation gives
\bea
d_{\Kbarbar}^2(\zeta_n,\mu) = \SE_{\Kbarbar}(\mu)-2\,\wb_n\TT \pb_{\Kbarbar, \mu}(\bX_n)+\wb_n\TT\Kbarbarb(\bX_n)\wb_n \,,
\eea
where $\pb_{\Kbarbar, \mu}(\bX_n)=\left[P_{\Kbarbar,\mu}(\bx^{(m+1)}),\ldots,P_{\Kbarbar,\mu}(\bx^{(m+n)})\right]\TT$, with $P_{\Kbarbar,\mu}(\bx)$ defined by \eqref{eq:potnu} in Appendix~A, $\{\Kbarbarb(\bX_n)\}_{i,j}=\Kbarbar(\bx^{(m+i)},\bx^{(m+j)})$ for any $i,j=1,\ldots,n$, and $\SE_{\Kbarbar}(\mu)=\int_{\X^2} \Kbarbar(\bx,\bx')\,\dd\mu(\bx)\dd\mu(\bx')$.

When $\bX_m\cap\bX_n=\emptyset$, the $n\times n$ matrix $\Kb_{|m}(\bX_n)$, whose element $i,j$ equals $K_{|m}(\bx^{(m+i)},\bx^{(m+j)})$, is positive definite. The elementwise (Hadamard) product $\Kb_{|m}(\bX_n)\circ\Kb_{|m}(\bX_n)$ is thus positive definite too, implying that $\Kbarbarb(\bX_n)$ is positive definite. The optimal weights $\wb_n^*$ minimizing
$d_{\Kbarbar}^2(\zeta_n,\mu)$ are thus 
\be\label{Eq:weights}
\wb_n^* = \Kbarbarb^{-1}(\bX_n)\pb_{\Kbarbar, \mu}(\bX_n) \,.
\ee
We shall denote by $\zeta_n^*$ the measure supported on $\bX_n$ with the optimal weights \eqref{Eq:weights} and
\be
Q_{n*}^2 &=& 1 - \frac{\ISE_{\zeta_n^*}(\bX_m,\by_m)}
{\frac1n \sum_{i=1}^n  \left[y(\bx^{(m+i)})-\overline{y}_n\right]^2} \nonumber \\ &=&  1 - \frac{\sum_{i=1}^n w_i^* \left[ y(\bx^{(m+i)})-\eta_m(\bx^{(m+i)})\right]^2}
{\frac1n \sum_{i=1}^n  \left[y(\bx^{(m+i)})-\overline{y}_n\right]^2} \,,
\label{Q2-est*}
\ee
with $\overline{y}_n=(1/n)\,\sum_{i=1}^n y(\bx^{(m+i)})$.
Notice that the weights $w_i^*$ do not depend on the variance parameter $\sigma^2$ of the GP model. 

\begin{remark}
When $\bX_n$ is constructed by kernel herding, see Section~\ref{S:KH}, $K$ can be chosen identical to the kernel used there. This will be the case in Sections~\ref{sec:4} and \ref{sec:5}, but it is not mandatory. 

Conversely, one may think of choosing a design $\bX_n$ that minimizes $d_{\Kbarbar}^2(\xi_n,\mu)$, or $d_{\Kbarbar}^2(\zeta_n^*,\mu)$, also with the objective to obtain a precise estimation of $\ISE_\mu(\bX_m,\by_m)$. This can be achieved by kernel herding using the kernel $\Kbarbar$ and is addressed in \cite{PR2021a}. However, the numerical results presented there show that the precise choice of the test set $\bX_n$ has a marginal effect compared to the effect of non-uniform weighting with $\wb_n^*$, provided that $\bX_n$ fills the holes left in $\X$ by the training design $\bX_m$.
\fin
\end{remark}

\begin{remark} 
When the observations $y(\bx^{(i)})$, $i=1,\ldots,n$, are available at the validation stage, an alternative version of $\widehat Q^2_n$ would be
\be\label{eq:Q2testprime}
\widehat Q'^2_n = 1 - \frac{ \sum_{\bx \in \bX_n}  \left[ y(\bx)-\eta_m(\bx)\right]^2}{\sum_{\bx \in \bX_n}  \left[y(\bx)-\overline{y}_m\right]^2}\,,
\ee
where $\overline{y}_m=(1/m)\,\sum_{i=1}^m y(\bx^{i})$, which compares the performance on the test set of two predictors $\eta_m$ and $\overline{y}_m$ based on the same training set. It is then possible to also apply a weighting procedure to the {\em denominator} of $\widehat Q'^2_n$,
\bea
D_{\xi_n}(\bX_m,\by_m) = \frac1n\,\sum_{i=1}^n  \left[y(\bx^{(m+i)})-\overline{y}_m\right]^2\ ,
\eea
in order to make it resemble its idealized version $V'_\mu(\by_m) = \int_\X \left[y(\bx)- \overline{y}_m \right]^2 \, \dd\mu(\bx)$. 
The GP model is now
$\mve_m(\bx)=y(\bx)-\overline{y}_m\sim\GP(\eta_m(\bx)-\overline{y}_m,\sigma^2\,K_{|m})$. Similar developments to those used above for the numerator of $\widehat Q^2_n$ yield
\bea
\overline{\Delta'}^2(\xi_n,\mu;\bX_m,\by_m) &=& \EE \left\{ \left[ D_{\xi_n}(\bX_m,\by_m) - V_\mu(\by_m) \right]^2 \right\}\\
&=& \EE\left\{ \int_{\X^2}  \mve_m^2(\bx)\mve_m^2(\bx')\,\dd(\xi_n-\mu)(\bx)\dd(\xi_n-\mu)(\bx')\right\} \\
&=& \sigma^2\,d_{\Kbarbarp}^2(\xi_n,\mu)\,,
\eea
where 
\bea
    \Kbarbarp(\bx, \bx') &=& \Kbarbar(\bx, \bx') + [\eta_m(\bx)-\overline{\by}_m]^2[\eta_m(\bx')-\overline{\by}_m]^2 \\
    &&+\,[\eta_m(\bx)-\overline{\by}_m]^2K_{|m}(\bx',\bx') + [\eta_m(\bx')-\overline{\by}_m]^2K_{|m}(\bx,\bx) \\ && +\,4\,[\eta_m(\bx)-\overline{\by}_m][\eta_m(\bx')-\overline{\by}_m]K_{|m}(\bx,\bx') \,.
\eea
We can then substitute $D_{{\zeta_n'}^*}(\bX_m,\by_m)$ for $D_{\xi_n}(\bX_m,\by_m)$ in \eqref{eq:Q2testprime}, where ${\zeta_n'}^*$ allocates the weights
$
{\wb_n'}^* = \Kbarbarpb^{-1}(\bX_n)\pb_{\Kbarbarp, \mu}(\bX_n) 
$
to the $n$ points in $\bX_n$, 
with $\{\Kbarbarpb(\bX_n)\}_{i,j}=\Kbarbarp(\bx^{(m+i)},\bx^{(m+j)})$, $i,j=1,\ldots,n$. 
\fin
\end{remark}

\section{Test-set construction}
\label{sec:2}

In the previous section we assumed  the test set as given, and proposed a method to estimate $\ISE_\mu(\bX_m,\by_m)$ by a weighted sum of the residuals. In this section we address the choice of the test set.

Below we give an overview of the three methods used in this paper, all relying on the concept of space-filling design \cite{fanli06, promul12}. While most methods for the construction of such designs choose all points simultaneously, the methods we consider are incremental, selecting one point at a time.
 
Our objective is to construct an ordered test set of size $n$, denoted by $\bX_n=\left\{\bx^{(1)},\ldots,\bx^{(n)}\right\}\subset\X$. When there is no restriction on the choice of $\bX_n$, the advantage of using an incremental construction is that it can be stopped once the estimation of the predictivity of an initial model, built with some given design $\bX_m$, is considered sufficiently accurate. In case the conclusion is that model predictions are not reliable enough, the full design $\bX_{m+n}=\bX_m\cup\bX_n$ and the associated observations $\by_{m+n}$ can be used to update the model. 
This updated model can then be tested at additional design points, elements of a new test set to be constructed. All methods presented in this section (except the Fully Sequential Space-Filling method) are implemented in the Python package \pyvar{otkerneldesign} \footnote{\url{https://pypi.org/project/otkerneldesign/}} which is based on the OpenTURNS library for uncertainty quantification \cite{baudut17}.

\subsection{Fully-Sequential Space-Filling design}\label{S:FSSF}

The Fully-Sequential Space-Filling forward-reflected (FSSF-fr) algorithm \cite{shaapl21} relies on the CADEX algorithm \cite{kensto69} (also called the ``coffee-house'' method~\cite{mul07}). It constructs a sequence of nested designs in a bounded set $\X$ by sequentially selecting a new point $\bx$ as far away as possible from the $\bx^{(i)}$ previously selected. New inserted points are selected within a set of candidates $\S$ which may coincide with $\X$ or be a finite subset of $\X$ (which simplifies the implementation, only this case is considered here). 
The improvement of FSSF-fr when compared to CADEX is that new points are selected {\em at the same time} far from the previous design points as well as  far from the boundary of $\X$.  

The algorithm is as follows:
\begin{enumerate}
    \item Choose $\S$, a finite set of candidate  points in $\X$, with size $N \gg n$ in order to allow a fairly dense covering of $\X$. When $\X=[0,1]^d$, \cite{shaapl21} recommends to take $\S$ equal to the first $N=1\,000\,d+2\,n$ points of a Sobol sequence in $\X$. 
    
    \item Choose the first point $\bx^{(1)}$ randomly in $\S$ and define $\bX_1 = \{\bx^{(1)}\}$. 
    
    \item At iteration $i$, with $\bX_i = \{\bx^{(1)},\ldots, \bx^{(i)}\}$, select
    \begin{equation}\label{eq:FSSF-fr}
        \bx^{(i+1)} \in \Argmax_{\bx \in \S \setminus \bX_i} \left[ \min\left( \min_{j \in \{1,\ldots,i\}} \|\bx-\bx^{(j)}\|, \sqrt{2}\,d \,\mbox{dist}(\bx,R(\bx)) \right) \right]\,, 
    \end{equation}
    where $R(\bx)$ is the symmetric of $\bx$ with respect to its nearest boundary of $\X$,
    and set $\bX_{i+1} = \bX_i \bigcup \bx^{(i+1)}$.
    
    \item Stop the algorithm when $\bX_n$ has the required size.
\end{enumerate}

The standard coffee-house (greedy packing) algorithm simply uses $\bx^{(i+1)} \in \Argmax_{\bx \in \S \setminus \bX_i}  \min_{j \in \{1,\ldots,i\}} \|\bx-\bx^{(j)}\|$. The role of the reflected point $R(\bx)$ is to avoid selecting $\bx^{(i+1)}$ too close to the boundary of $\X$, which is a major problem with standard coffee-house, especially when $\X=[0,1]^d$ with $d$ large. 
The factor $\sqrt{2}d$ in \eqref{eq:FSSF-fr} proposed in \cite{shaapl21} sets a balance between distance to the design $\bX_i$ and distance to the boundary of $\X$. Another scaling factor, depending on the target design  size $n$ is proposed in \cite{nogpro21}.
 
FSSF-fr is entirely based on geometric considerations and implicitly assumes that the selected set of points should cover $\X$ evenly. However, in the context of uncertainty quantification \cite{smi14} it frequently happens that the distribution $\mu$ of the model inputs is not uniform. It is then desirable to select a test set representative of  $\mu$. This can be achieved through the inverse probability integral transform: FSSF-fr constructs $\bX_n$ in the unit hypercube $[0,1]^d$, and an  ``isoprobabilistic'' transform $T:[0,1]^d \rightarrow \X$ is then applied to the points in $\bX_i$, $T$ being such that, if $\boldsymbol{U}$ is a random variable  uniform on $[0,1]^d$, then $T(\boldsymbol{U})$ follows the target distribution $\mu$. The transformation can be applied to each input separately when $\mu$ is the product of its marginals, a situation considered in our second test-case of Section~\ref{sec:4}, but is more complicated in other cases, see \cite[Chap.~4]{lemcha09}. 
Note that FSSF-fr operates in the bounded set $[0,1]^d$ even if the support of $\mu$ is unbounded. 
The other two algorithms presented in this section are able to directly choose points representative of a given distribution $\mu$ and do not need to resort to such a transformation.

\subsection{Support points}\label{S:SP}

Support points \cite{makjos18} are such that their associated empirical distribution $\xi_n$ has minimum Maximum-Mean-Discrepancy (MMD) with respect to $\mu$ for the energy-distance kernel of Sz\'ekely and Rizzo~\cite{szeriz04,szeriz13},  
\begin{equation}\label{eq:kE}
K_E(\bx,\bx') = \frac12\, \left(\Vert \bx \Vert + \Vert \bx' \Vert - \Vert \bx-\bx' \Vert\right)\,.
\end{equation} 
The squared MMD between $\xi_n$ and $\mu$ for the distance kernel equals
\begin{equation}\label{eq:energySP}
d_{K_E}^2(\xi_n,\mu) = \frac{2}{n} \sum_{i=1}^{n} \EE \|\bx^{(i)}-\zeta\| - \frac{1}{n^2} \sum_{i=1}^{n}\sum_{j=1}^{n}  \|\bx^{(i)}-\bx^{(j)}\| -  \EE \|\zeta-\zeta'\| \,,
\end{equation}
where $\zeta$ and $\zeta'$ are independently distributed with $\mu$; see \cite{sejsri13}.
A key property of the energy-distance kernel is that it is characteristic \cite{srigre10}: for any two probability distributions $\mu$ and $\xi$ on $\X$, $d_{K_E}^2(\mu,\xi)$ equals zero if and only if $\mu=\xi$, and so it defines a norm in the space of probability distributions.
Compared to more heuristic methods for solving quantization problems, support points
benefit from the theoretical guarantees of MMD minimization in terms of convergence of $\xi_n$ to $\mu$ as $n\to\infty$. 

As $\EE \|\bx^{(i)}-\zeta\|$ is not known explicitly, in practice $\mu$ is replaced by its empirical version $\mu_N$ for a given large-size sample $(\bx'^{(k)})_{k=1\ldots N}$. 
The support points $\bX_n^s$ are then given by
\begin{equation}\label{eq:SPestim}
\bX_n^s \in \Argmin_{\bx^{(1)},\ldots,\bx^{(n)}} \left( \frac{2}{n N} \sum_{i=1}^{n} \sum_{k=1}^N \|\bx^{(i)}-\bx'^{(k)}\| - \frac{1}{{n}^2} \sum_{i=1}^{n}\sum_{j=1}^{n} \|\bx^{(i)}-\bx^{(j)}\| \right) \,.
\end{equation}
The function to be minimized can be written as a difference of functions convex in $\bX_n$, which yields a difference-of-convex program. 
In \cite{makjos18}, a majorization-minimization procedure, efficiently combined with resampling, is applied to the construction of large designs (up to $n=10^4$) in high dimensional spaces (up to $d=500$). The examples treated clearly show that support points are distributed in a way that matches $\mu$ more closely than Monte-Carlo and quasi-Monte Carlo samples \cite{fanli06}.

The method can be used to split a dataset into a training set and a test set \cite{josvak22}: the $N$ points $\bX_N$ in \eqref{eq:SPestim} are those from the dataset, $\bX_n^s$ gives the test set and the other $N-n$ points are used for training. There is a serious additional difficulty though, as choosing $\bX_n^s$ among the dataset corresponds to a difficult combinatorial optimization problem. A possible solution is to perform the optimization in a continuous domain $\X$ and then choose $\bX_n^s$ that corresponds to the closest points in $\bX_N$ (for the Euclidean distance) to the continuous solution obtained \cite{josvak22}. 

The direct determination of support points through \eqref{eq:SPestim} does not allow the construction of a nested sequence of test sets. One possibility would be to solve \eqref{eq:SPestim} sequentially, one point at a time, in a continuous domain, and then select the closest point within $\bX_N$ as the one to be included in the test set. 
We shall use a different approach here, based on the greedy minimization of the MMD \eqref{eq:energySP} for the candidate set $\S=\bX_N$: at iteration $i$, the algorithm chooses
\begin{equation}\label{eq:GreedySP}
\bx_{i+1}^s \in \Argmin_{\bx\in\S} \left( \frac{1}{N} \sum_{k=1}^N \|\bx-\bx'^{(k)}\| - \frac{1}{i+1} \sum_{j=1}^{i} \|\bx-\bx^{(j)}\| \right) \,.
\end{equation}
The method requires the computation of the $N(N-1)/2$ distances $\|\bx^{(i)}-\bx^{(j)}\|$, $i,j=1,\ldots,N$, $i\neq j$, which hinders its applicability to large-scale problems (a test-case with $N=1\,000$ is presented in Section~\ref{sec:5}). Note that we consider support points in the input space $\X$ only, with $\X\subseteq\R^d$, in contrast with \cite{josvak22} which considers couples $(\bx^{(i)},y(\bx^{(i)}))$ in $\R^{d+1}$ to split a given dataset into a training set and a test set. 

Greedy MMD minimization can be applied to other kernels than the distance kernel \eqref{eq:kE}, see \cite{TeymurGRO2021, pro21}. In the next section we consider the closely related method of Kernel Herding (KH) \cite{chewel10}, which corresponds to a conditional-gradient descent   in the space of probability measures supported on a candidate set $\S$; see, e.g., \cite{prozhi20} and the references therein.

\subsection{Kernel herding}\label{S:KH}

Let $K$ be a positive definite kernel on $\X\times\X$. 
At iteration $i$ of kernel herding, with $\xi_i=(1/i)\sum_{j=1}^i \delta_{\bx^{(i)}}$ the empirical measure for $\bX_i$, the next point $\bx_{i+1}$ minimizes the directional derivative $F_K(\xi_i,\mu,\delta_\bx)$ of the squared MMD $d_K^2(\xi,\mu)$ at $\xi=\xi_i$ in the direction of the delta measure $\delta_\bx$, see Appendix~A. Direct calculation gives
$F_K(\xi_i,\mu,\delta_\bx)=2\,[P_{K,\xi}(\bx)-P_{K,\mu}(\bx)]-2\,\int_\X [P_{K,\xi}(\bx)-P_{K,\mu}(\bx)]\,\dd\xi(\bx)$, with $P_{K,\xi}(\bx)$ (resp.\ $P_{K,\mu}(\bx)$) the potential of $\xi$ (resp.\ $\mu$) at $\bx$, see \eqref{eq:potnu}, and thus
\begin{eqnarray}\label{eq:KH}
\bx_{i+1} &\in& \Argmin_{\bx\in\S} \left[P_{K,\xi_i}(\bx)-P_{K,\mu}(\bx)\right] \,,
\end{eqnarray}
with $\S\subseteq\X$ a given candidate set. 
Here, $P_{K,\xi_i}(\bx) = (1/i)\, \sum_{j=1}^i K(\bx, \bx^{(j)})$. When an empirical measure $\mu_N$ based on a sample $(\bx'^{(k)})_{k=1\ldots N}$ is substituted for $\mu$, we get 
$P_{K,\mu_N}(\bx) = (1/N)\, \sum_{k=1}^N K(\bx, \bx'^{(k)})$, which gives
\begin{eqnarray*}
\bx_{i+1} &\in& \Argmin_{\bx\in\S} \left[ \frac{1}{i} \sum_{j=1}^{i} K(\bx,\bx^{(j)}) - \frac{1}{N} \sum_{k=1}^N K(\bx,\bx'^{(k)}) \right] \,.
\end{eqnarray*}
When $K$ is the energy-distance kernel \eqref{eq:kE} we thus obtain \eqref{eq:GreedySP} with a factor $1/i$ instead of $1/(i+1)$ in the second sum. 

The candidate set $\S$ in \eqref{eq:KH} is arbitrary and can be chosen as in Section~\ref{S:FSSF}. 
A neat advantage of kernel herding over support points in that the potential $P_{K,\mu}(\bx)$ is sometimes explicitly available. When $\S=\bX_N$, this avoids the need to calculate the $N(N-1)/2$ distances $\|\bx^{(i)}-\bx^{(j)}\|$ and thus allows application to very large sample sizes. This is the case in particular when $\X$ is the cross product of one-dimensional sets $\X_{[i]}$, $\X=\X_{[1]}\times\cdots\times\X_{[d]}$, $\mu$ is the product of 
its marginals
$\mu_{[i]}$ on the $\X_{[i]}$, $K$ is the product of one-dimensional kernels $K_{[i]}$, and the one-dimensional integral in  $P_{K_{[i]},\mu_{[i]}}(x)$ is known explicitly for each $i\in\{1,\ldots,d\}$. Indeed, for $\bx=(x_1,\ldots,x_d)\in\X$, we then have $P_{K,\mu}(\bx)=\prod_{i=1}^d P_{K_{[i]},\mu_{[i]}}(x_i)$; see \cite{prozhi20}. When $K$ is the product of Mat\'ern kernels with regularity parameter $5/2$ and correlation lengths $\theta_i$, $K(\bx,\bx') = \prod_{i=1}^{d} K_{5/2,\theta_{i}}(x_{i}-x'_i)$, with
\begin{equation}\label{eq:Matern5/2}
K_{5/2,\theta}(x-x')
=
\left(1 + \frac{\sqrt{5}}{\theta} |x - x'| + \frac{5}{3 \theta^2} (x - x')^2 \right)
\exp \left( - \frac{\sqrt{5}}{\theta} |x - x'| \right),
\end{equation}
the one-dimensional potentials are given in Appendix~B for $\mu_{[i]}$ uniform on $[0,1]$ or $\mu_{[i]}$ the standard normal $\SN(0,1)$.  
When no observation is available, which is the common situation at the design stage, the correlation lengths have to be set to heuristic values. We empirically found the values of the correlation lengths to have a large influence over the design. A reasonable choice for $\X=[0,1]^d$ is $\theta_i = n^{-1/d}$ for all $i$, with $n$ the target number of design points; see \cite{prozhi20}. 

\subsection{Numerical illustration}\label{S:numerical-1}

We apply FSSF-fr (denoted FSSF in the following), support points and kernel herding algorithms to the situation where a given initial design of size $m$ has to be completed by a series of additional points $\bx^{(m+1)},\ldots,\bx^{(m+n)}$. The objective is to obtain a full design $\bX_{m+n}$ that is a good quantization of a given distribution $\mu$. 

\begin{sidewaysfigure}
    \vspace{11cm}
    \begin{minipage}[t]{\textwidth}
        \includegraphics[width=0.32\textwidth]{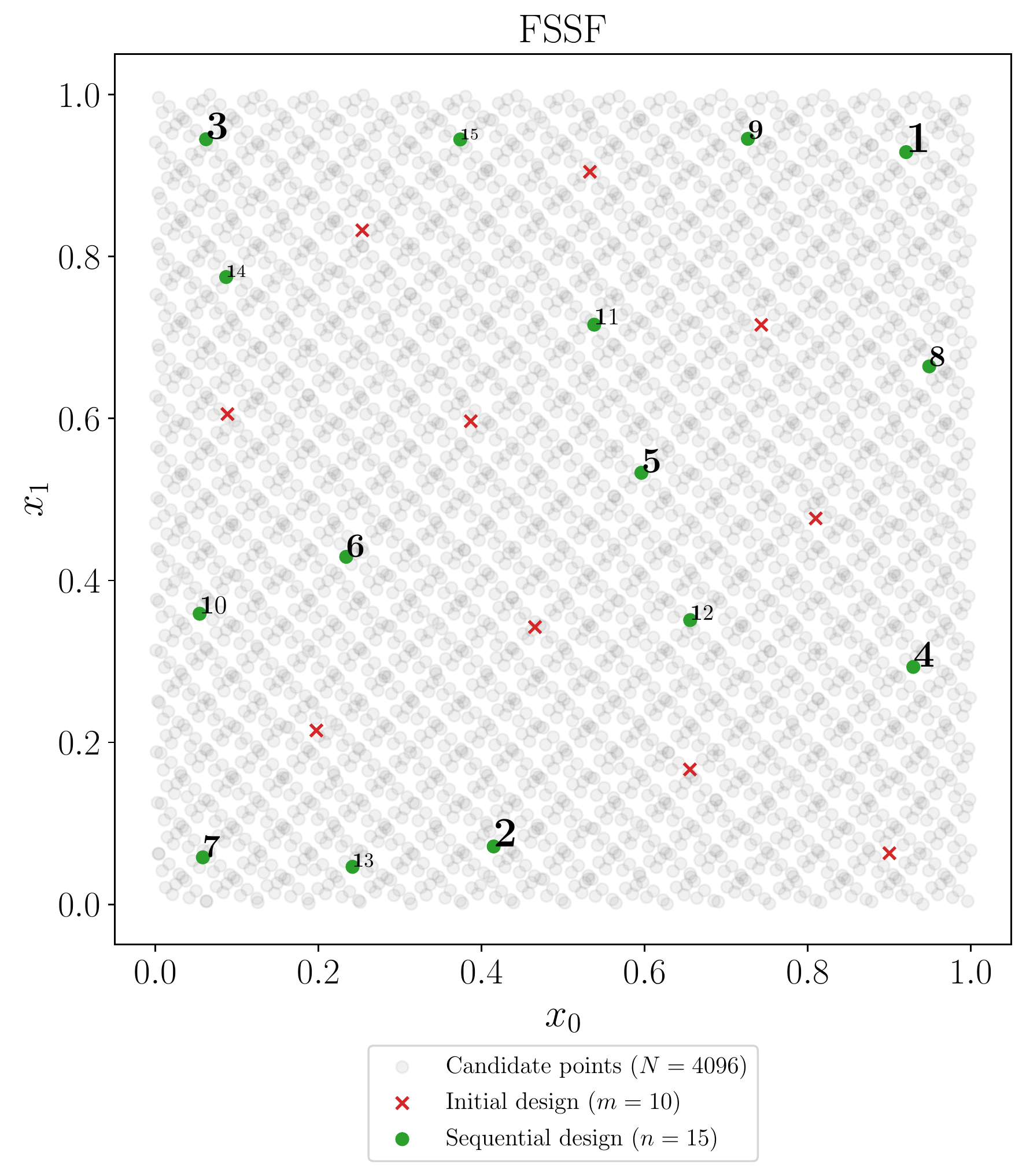}
        \quad
        \includegraphics[width=0.32\textwidth]{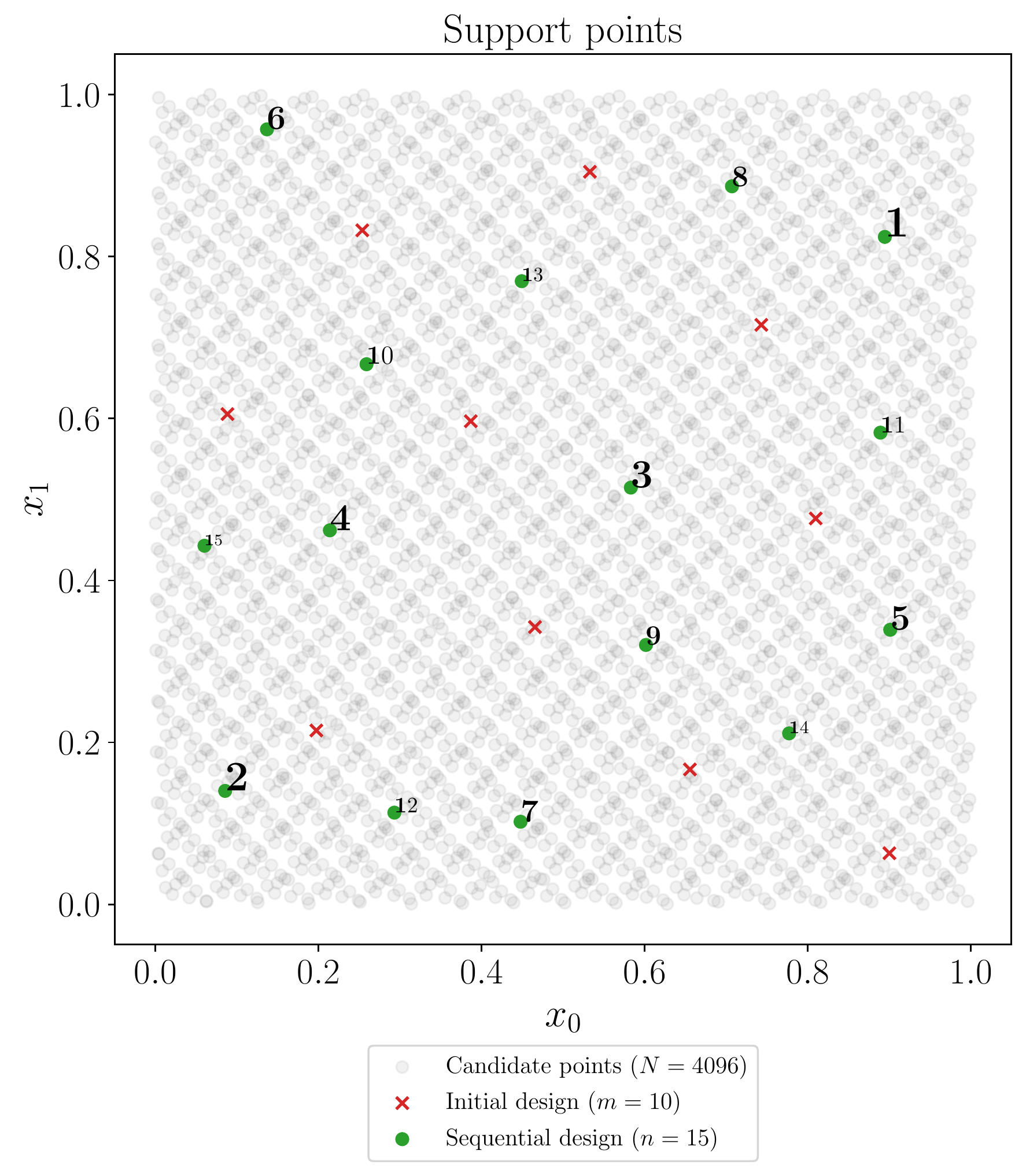}
        \quad
        \includegraphics[width=0.32\textwidth]{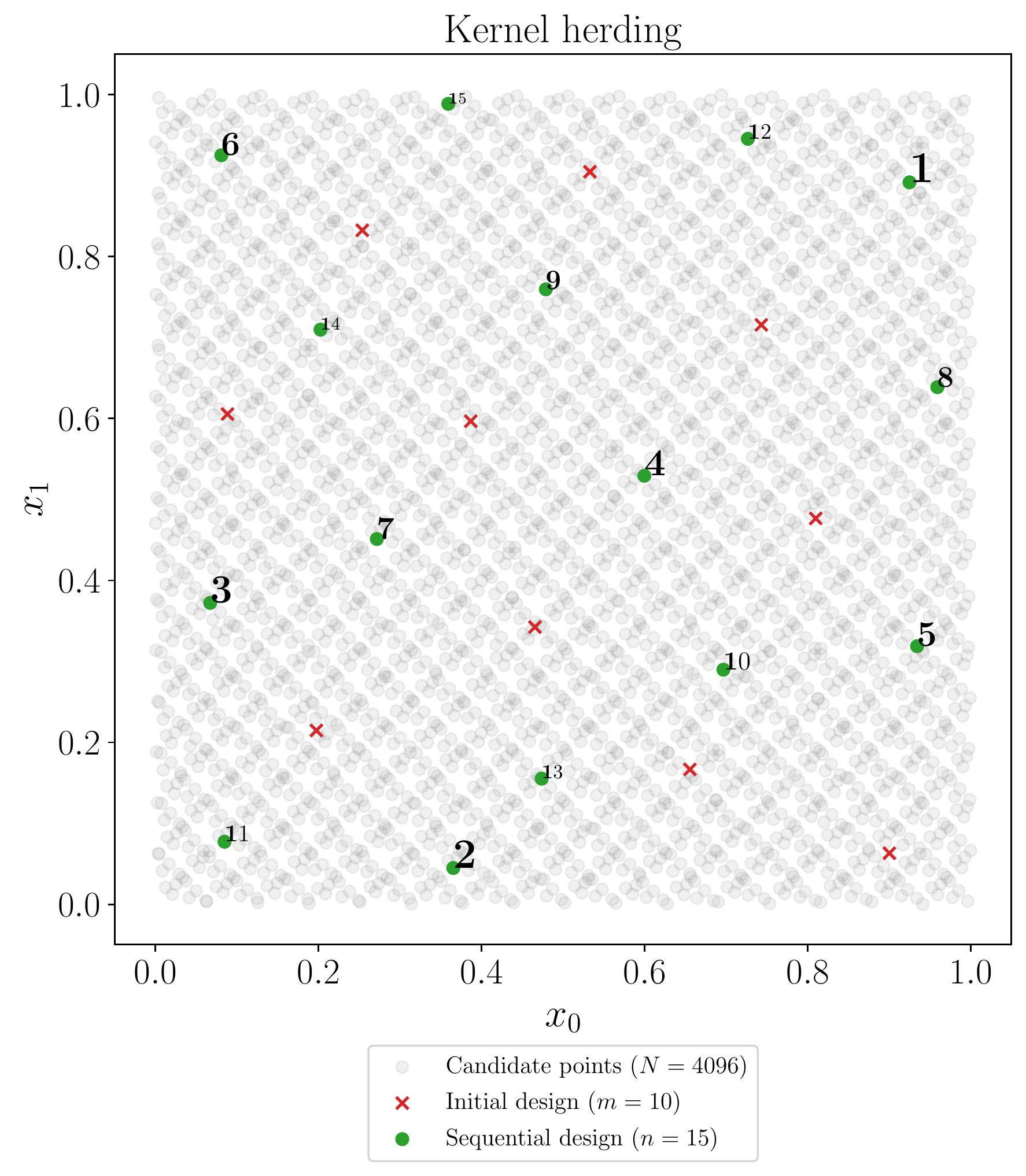}
        \caption{Additional points (ordered, green) complementing an initial design (red crosses), $\mu$ is uniform on $[0,1]$, the candidate points are in gray.}
        \label{fig:uniform_validation_designs}
    \end{minipage}
    \vspace{1cm}
     
    \begin{minipage}[t]{\textwidth}
        \includegraphics[width=0.32\textwidth]{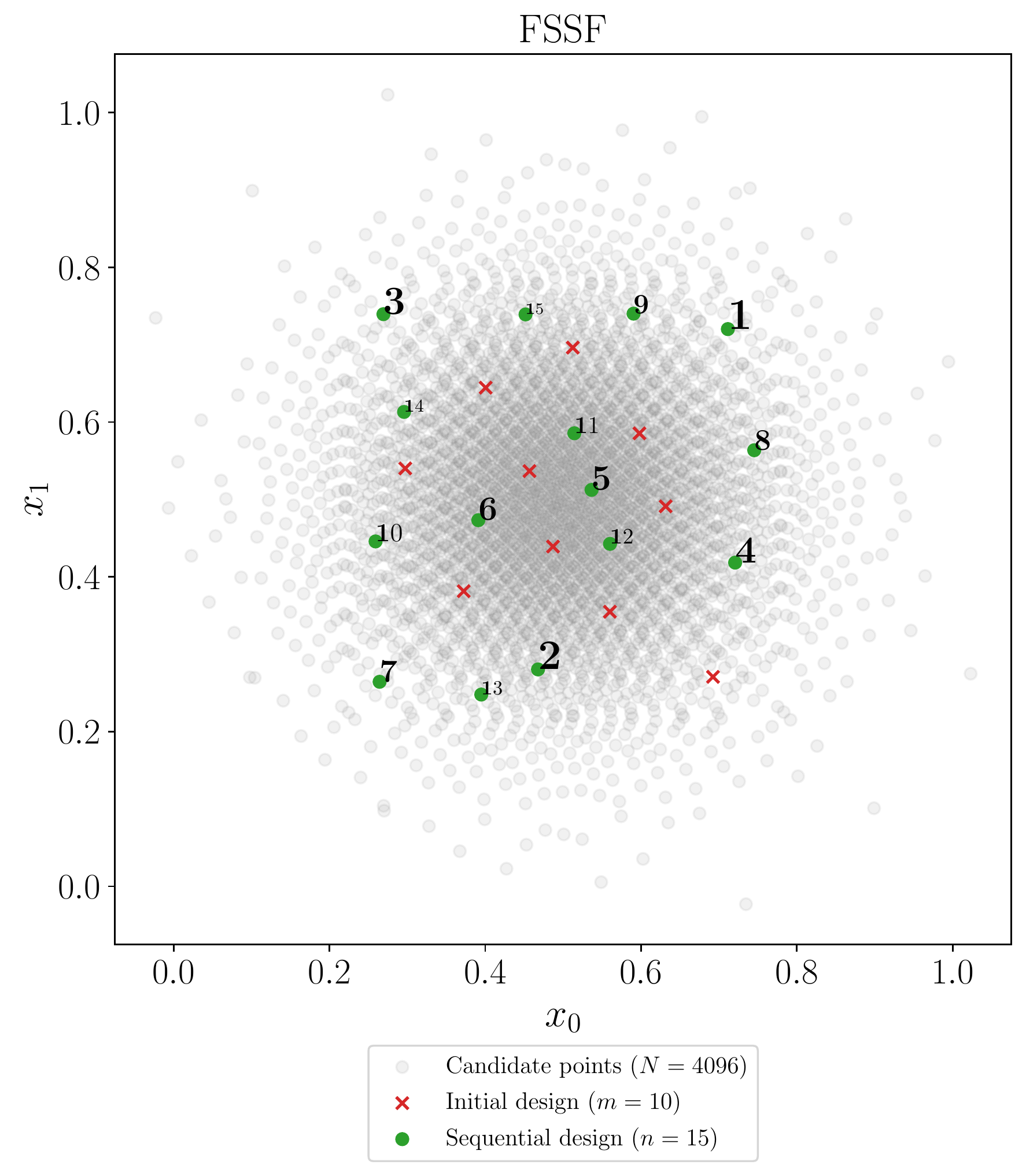}
        \quad
        \includegraphics[width=0.32\textwidth]{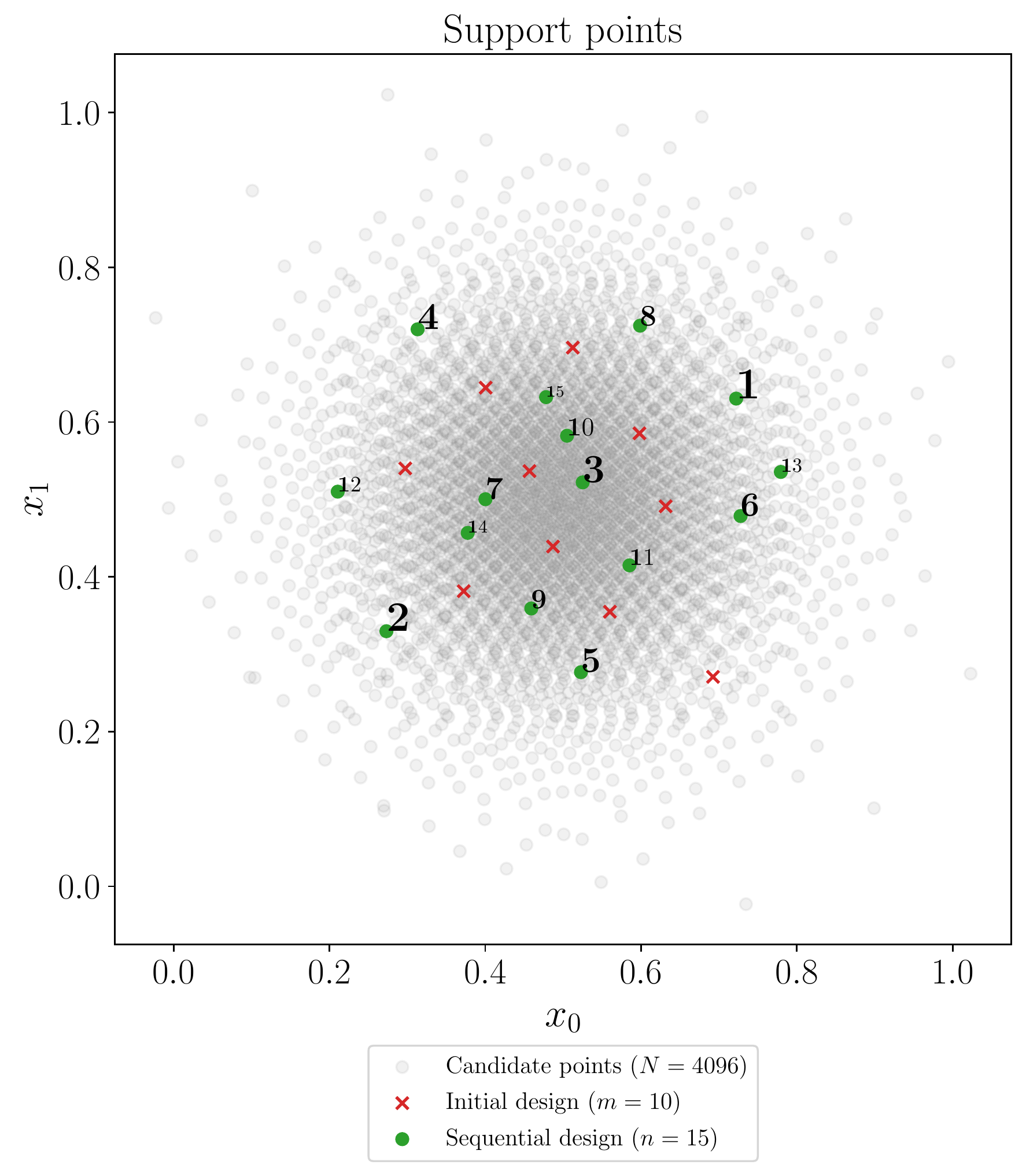}
        \quad
        \includegraphics[width=0.32\textwidth]{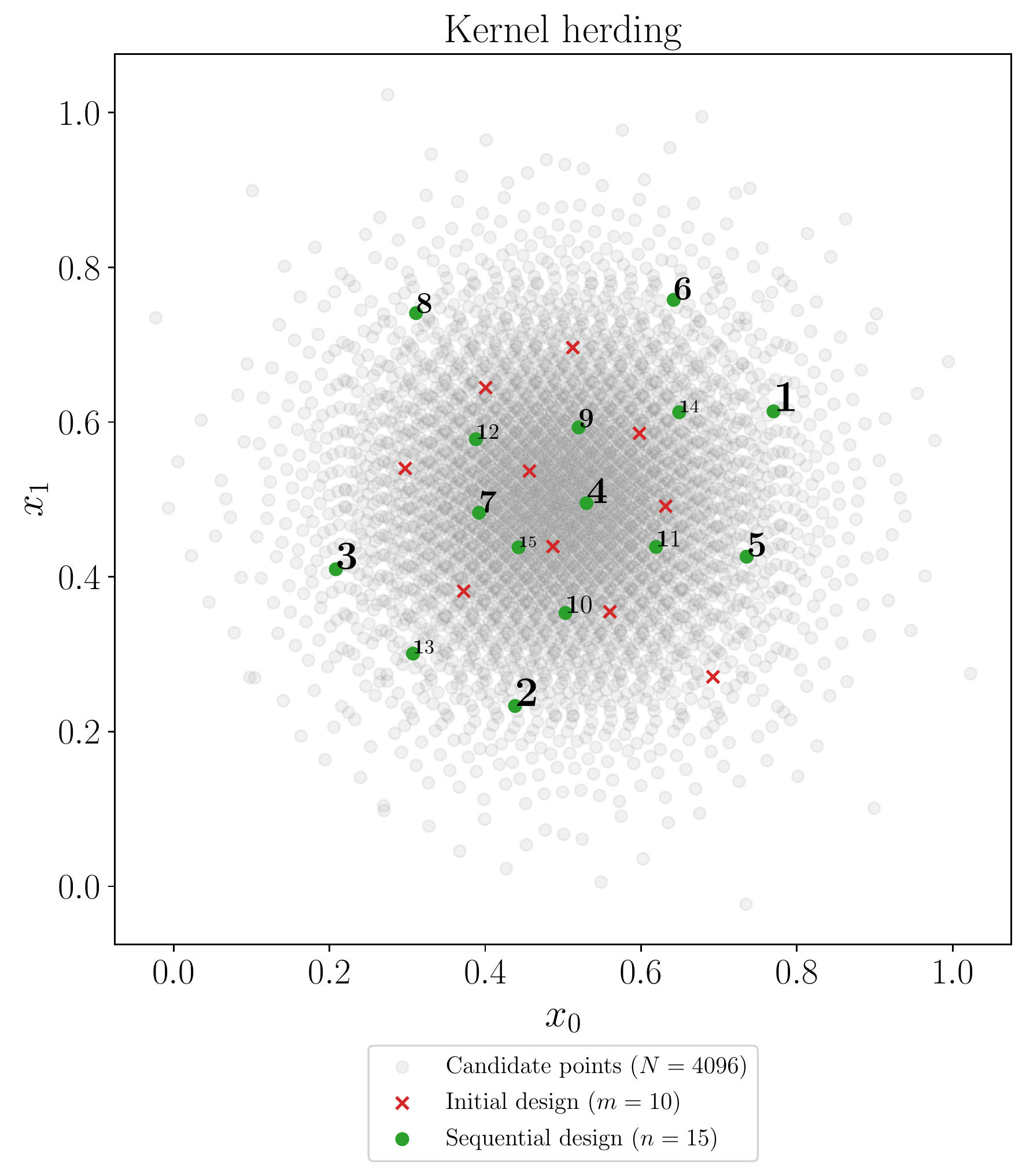}
        \caption{Additional points (ordered, green) complementing an initial design (red crosses), $\mu$ normal, the candidate points are in gray.}
        \label{fig:normal_validation_designs}
    \end{minipage}
\end{sidewaysfigure}

Figures~\ref{fig:uniform_validation_designs} and \ref{fig:normal_validation_designs}  correspond to $\mu$ uniform on $[0,1]^2$ and $\mu$ the standard normal distribution $\SN(\0b,\Ib_2)$, with $\Ib_2$ the 2-dimensional identity matrix, respectively. All methods  are applied to the same candidate set $\S$. 

The initial designs $\bX_m$ are chosen in the class of space-filling designs, well suited to initialize sequential learning strategies \cite{sanwil03}.
When $\mu$ is uniform, the initial design is a maximin Latin hypercube design \cite{mormit95} with $m=10$ and the candidate set is given by the $N=2^{12}$ first points $\Sb_N$ of a Sobol sequence in $[0,1]$. When $\mu$ is normal, the inverse probability transform method is first applied to $\Sb_N$ and $\bX_m$ (this does not raise any difficulty here as $\mu$ 
is the product of its marginals).
The candidate points $\S$ are marked in gray on Figures~\ref{fig:uniform_validation_designs} and \ref{fig:normal_validation_designs} and the initial design is indicated by the red crosses. The index $i$ of each  added test point $\bx^{(m+i)}$ is indicated (the font size decreases with $i$). In such a small dimension ($d=2$), a visual appreciation gives the impression that the three methods have comparable performance. 
We can notice, however, that FSSF tends to choose points closer to the boundary of $\S$ than the other two, and that  support points seem to sample more freely the holes of $\bX_m$ than kernel herding, which seems to be closer to a space-filling continuation of the training set.
We will come back to these designs when analysing the  quality of the resulting predictivity metric estimators in the next section.



\section{Numerical results I: construction of a training set and a test set}
\label{sec:4}

This section presents numerical results obtained on three different test-cases, in dimension 2 (test-cases 1 and 2) and 8 (test-case 3), for which $y(\bx)=f(\bx)$ with $f(\bx)$ having an easy to evaluate analytical expression, see Section \ref{sec:testCases}. This allows a good  estimation of $Q_{\mathrm{ideal}}^2(\mu)$ by $Q_{MC}^2=Q_{\mathrm{ideal}}^2(\mu_M)$, see \eqref{eq:Q2th}, where $\mu_M$ is the empirical measure for a large Monte-Carlo sample ($M=10^6$), that will serve as  reference when assessing the performance of each of the other estimators. We consider the validation designs built by FSSF, support points and kernel herding, presented in Sections~\ref{S:FSSF}, \ref{S:SP}, and \ref{S:KH}, respectively, and, 
for each one, we compare the performances obtained for both the uniform and the weighted estimator of Section~\ref{S:weighting}. 

\subsection{Test-cases} \label{sec:testCases}

The training  design $\bX_m$ and the set $\S$  of potential test set points are as in Section~\ref{S:numerical-1}. For test-cases 1 and 3, $\mu$ is the uniform measure on $\X=[0,1]^d$, with $d=2$ and $d=8$, respectively; $\bX_m$ is a maximin Latin hypercube design in $\X$, and $\S$ corresponds to the first $N$ points $\Sb_N$ of Sobol' sequence in $\X$, complemented by the $2^d$ vertices. In the second test-case, $d=2$, $\mu$ is  the normal distribution $\SN(\0b,\Ib_2)$, and the sets $\bX_m$ and $\Sb_N$ must be transformed as explained in section \ref{S:FSSF}. There are $N=2^{14}$ candidate points for test-cases 1 and 2 and $N=2^{15}$ for test-case 3 (this value is rather moderate for a problem in dimension 8, but using a larger $N$ yields numerical difficulties for support points; see Section~\ref{S:SP}). 

For each test-case, a GP regression model is fitted to the $m$ observations using ordinary Kriging \cite{raswil06} (a GP model with constant mean), with an 
anisotropic Matérn kernel with regularity parameter $5/2$: we substitute $[(\bx-\bx')\TT\Db(\bx-\bx')]^{1/2}$ for $|x-x'|$ in \eqref{eq:Matern5/2}, 
with $\Db$ a diagonal matrix with diagonal elements $1/\theta_i^2$, and the correlation lengths $\theta_i$ are estimated by maximum likelihood via a truncated Newton algorithm. 
All calculations were done using the Python package OpenTURNS for uncertainty quantification \cite{baudut17}. The kernel used for kernel herding is different and corresponds to the tensor product of one-dimensional Matérn kernels \eqref{eq:Matern5/2}, so that the potentials $P_{K,\mu}(\cdot)$ are known explicitly (see Appendix~B); the correlations lengths are set to $\theta=0.2$ in test-cases 1 and 3 ($d=2$) and to $\theta=0.7$ in test-case 3 ($d=8$).

Assuming that a model is classified, in terms of the estimated value of its predictivity index $Q^2$ as ``poor fitting'' if $Q^2\in[0.6, 0.8]$, ``reasonably good fitting'', when $Q^2\in(0.8,0.9]$, and ``very good fitting'' if $Q^2>0.9$, we selected, for each test-case three different sizes $m$ of the training set such that the corresponding models cover all three possible situations. 
For all test-cases, the impact of the size $n$ of the test set is studied in the range $n\in\{4 ,\ldots,50\}$.

\paragraph{Test-case 1.}

\begin{figure}
    \centering
    \includegraphics[width=0.49\textwidth]{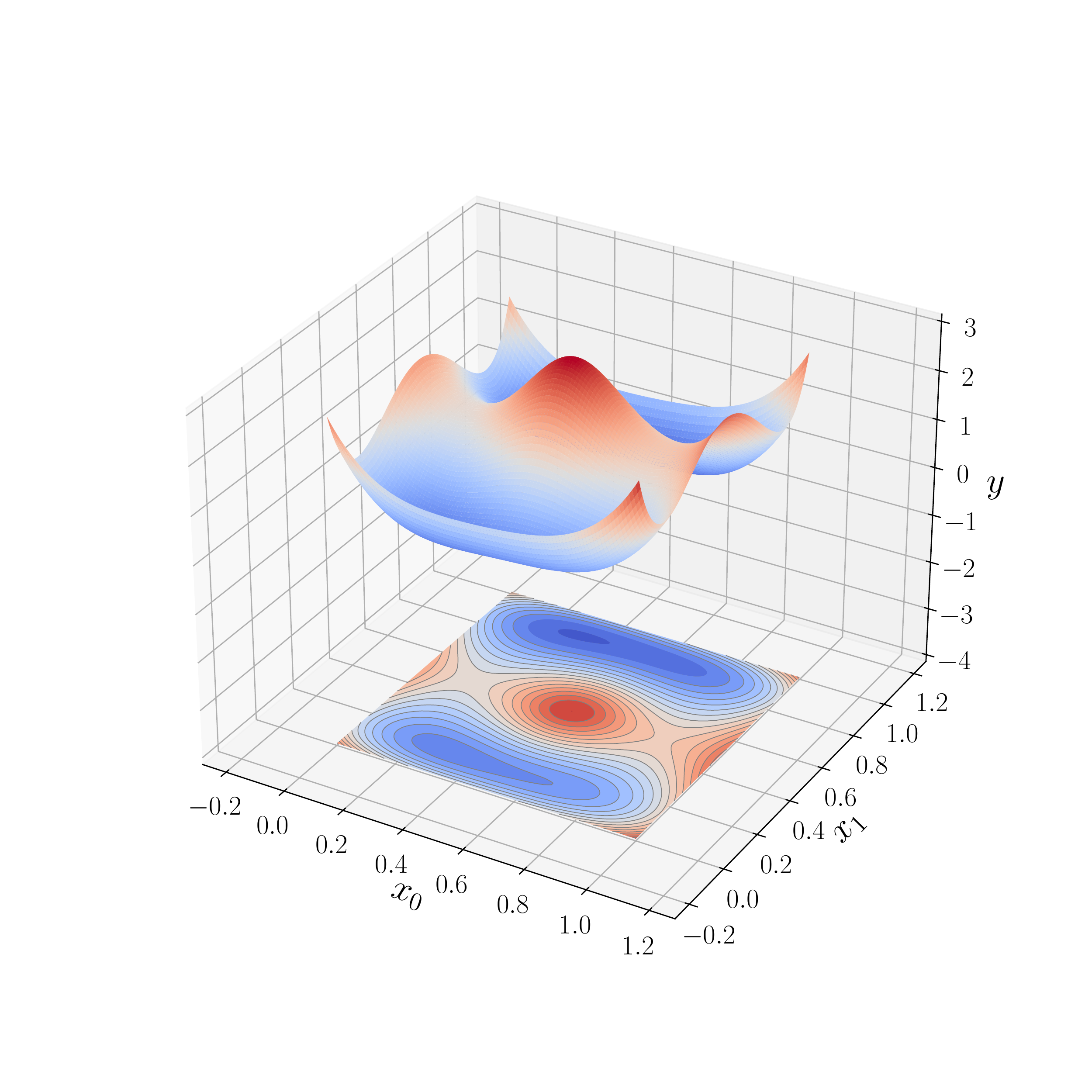}
        \includegraphics[width=0.49\textwidth]{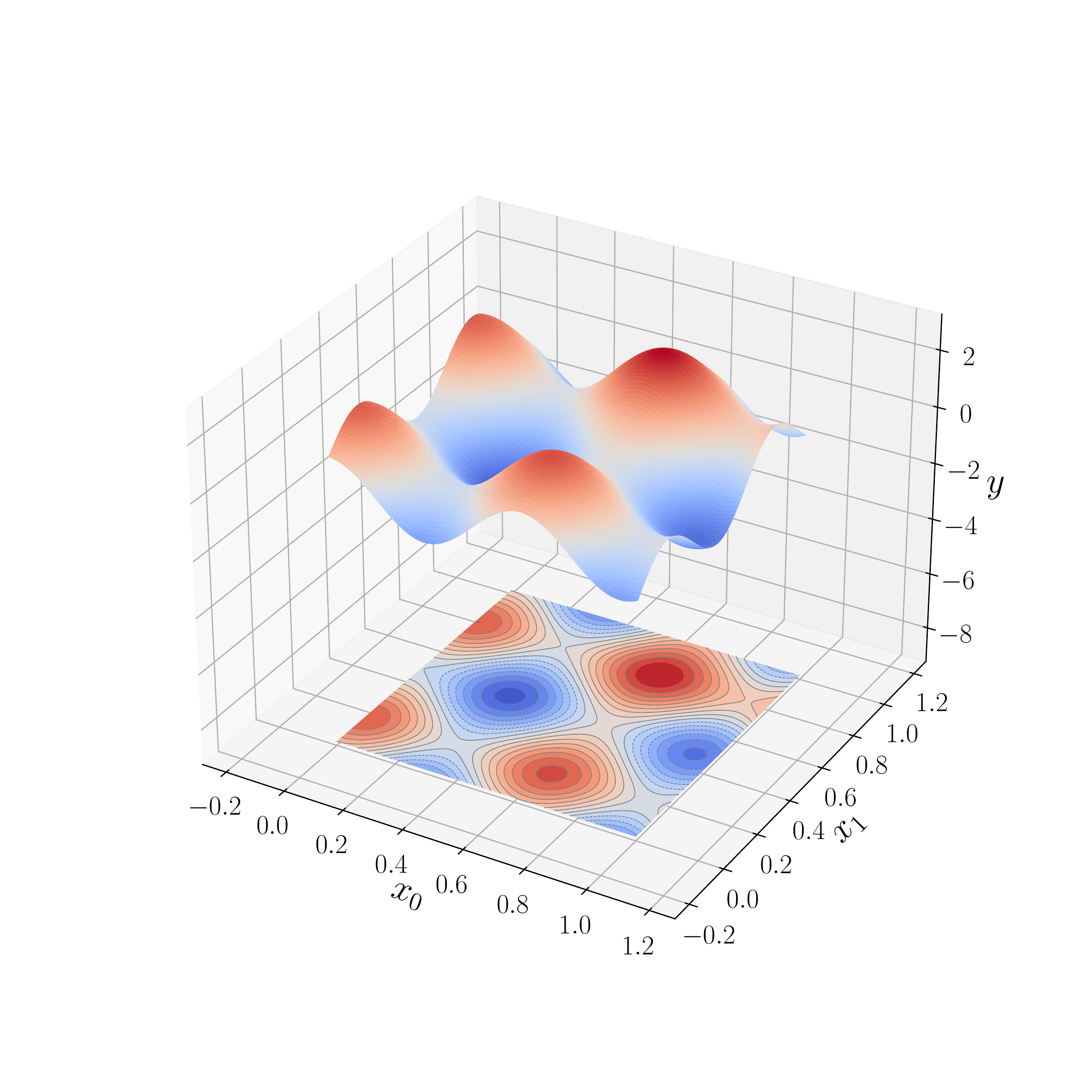}
    \caption{Left: $f_1(\bx)$ (test-case 1); right: $f_2(\bx)$ (test-case 2); $\bx\in\X=[0,1]^2$.}  
    \label{fig:f1&f2}
\end{figure}

This test function is $f_1(\bx) = h(2\,x_1 - 1, 2\,x_2 - 1)$, $(x_1,x_2) \in \X=[0,1]^2$, with
\begin{align*}
   h(u_1, u_2) =& \frac{\exp(u_1)}{5} - \frac{u_2}{5} + \frac{u_2^6}{3} + 4 u_2^4 - 4 u_2^2 + \frac{7u_1^2}{10} + u_1^4 + \frac{3}{4 u_1^2 + 4 u_2^2 + 1}\,.  
\end{align*}
Color coded 3d and contour plots of $f_1$ for $\bX\in\X$ are shown on the left panel of Figure~\ref{fig:f1&f2}, showing that 
the function is rather smooth, even if its behaviour along the boundaries of $\X$, in particular close to the vertices, may present difficulties for some regression methods. The size of the training set for this function are: $m\in\{5, 15, 30\}$.

\vspace{-0.6cm}
\vspace{-0.8cm}


\paragraph{Test-case 2.}

The second test function,  plotted in the right panel of Figure~\ref{fig:f1&f2} for $\bx\in[0,1]^2$, is 
\begin{align*}
     f_2(\bx) 
     &= \cos \left(5 + \frac{3}{2}x_1 \right) +  \sin \left(5 + \frac{3}{2}x_1 \right) 
     + \frac{1}{100} \left(5 + \frac{3}{2}x_1 \right) \left(5 + \frac{3}{2}x_2 \right)\,.
\end{align*}
Training set sizes for this test-case  are $m\in\{8, 15, 30\}$.

\paragraph{Test-case 3.}

The third function is the so-called ``gSobol'' function, defined over $\X=[0,1]^8$ by
\begin{equation*}
    f_3(\bx) = \prod_{i=1}^{8} \frac{|4x_i - 2| + a_i}{1 + a_i}, \qquad a_i = i^2 \,.
\end{equation*}
This parametric function is very versatile as both the dimension of its input space and the coefficients $a_i$ can be freely chosen. The sensitivity to input variables is determined by the $a_i$: the larger $a_i$ is, the less $f$ is sensitive to $x_i$. Larger  training sets are considered for this test-case: $m\in\{15, 30, 100\}$.

\subsection{Results and analysis}

\begin{sidewaysfigure}
    \vspace{11cm}
    
    \begin{minipage}[t]{\textwidth}
        \includegraphics[width=0.3\textwidth]{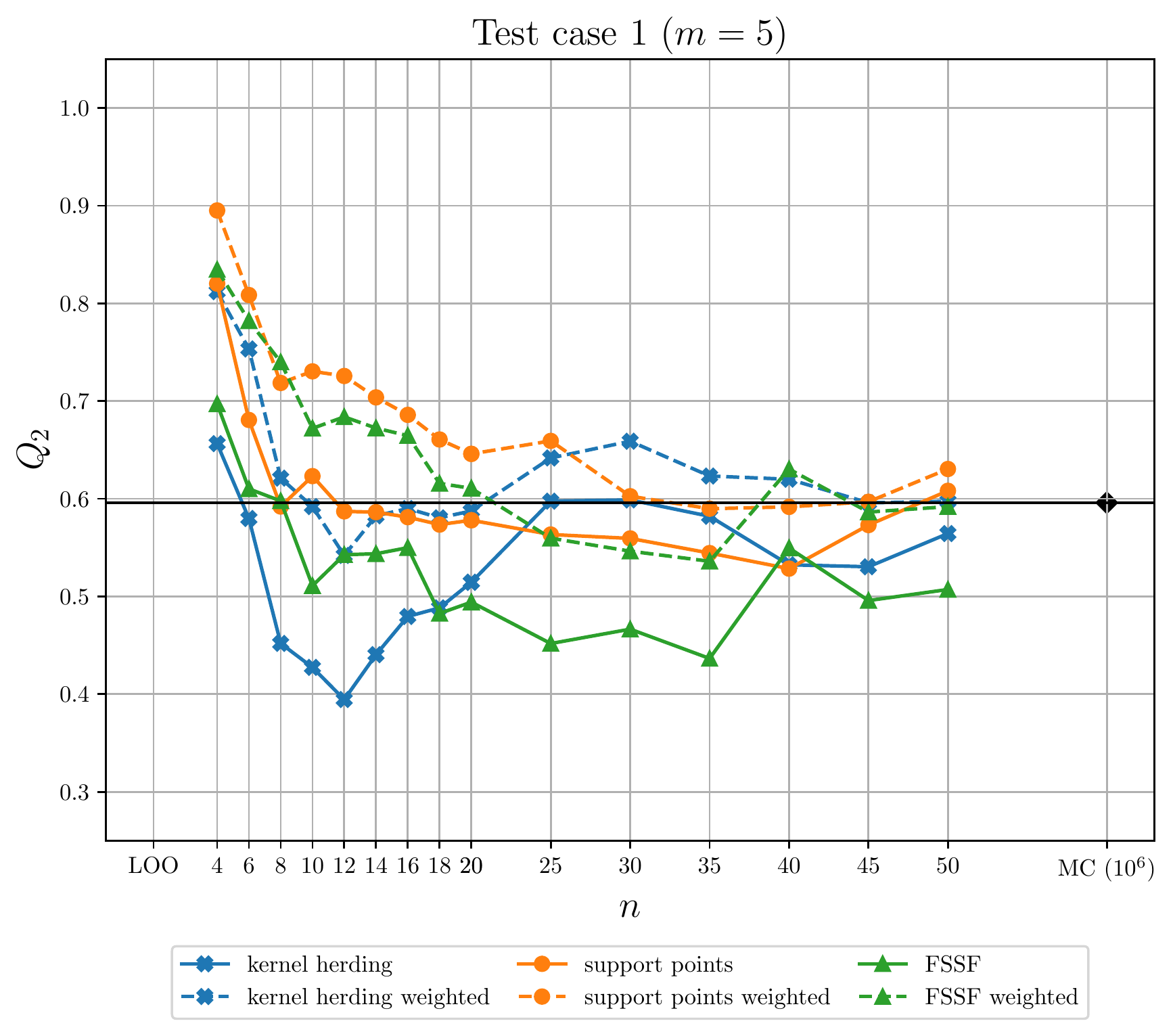}
        \quad
        \includegraphics[width=0.3\textwidth]{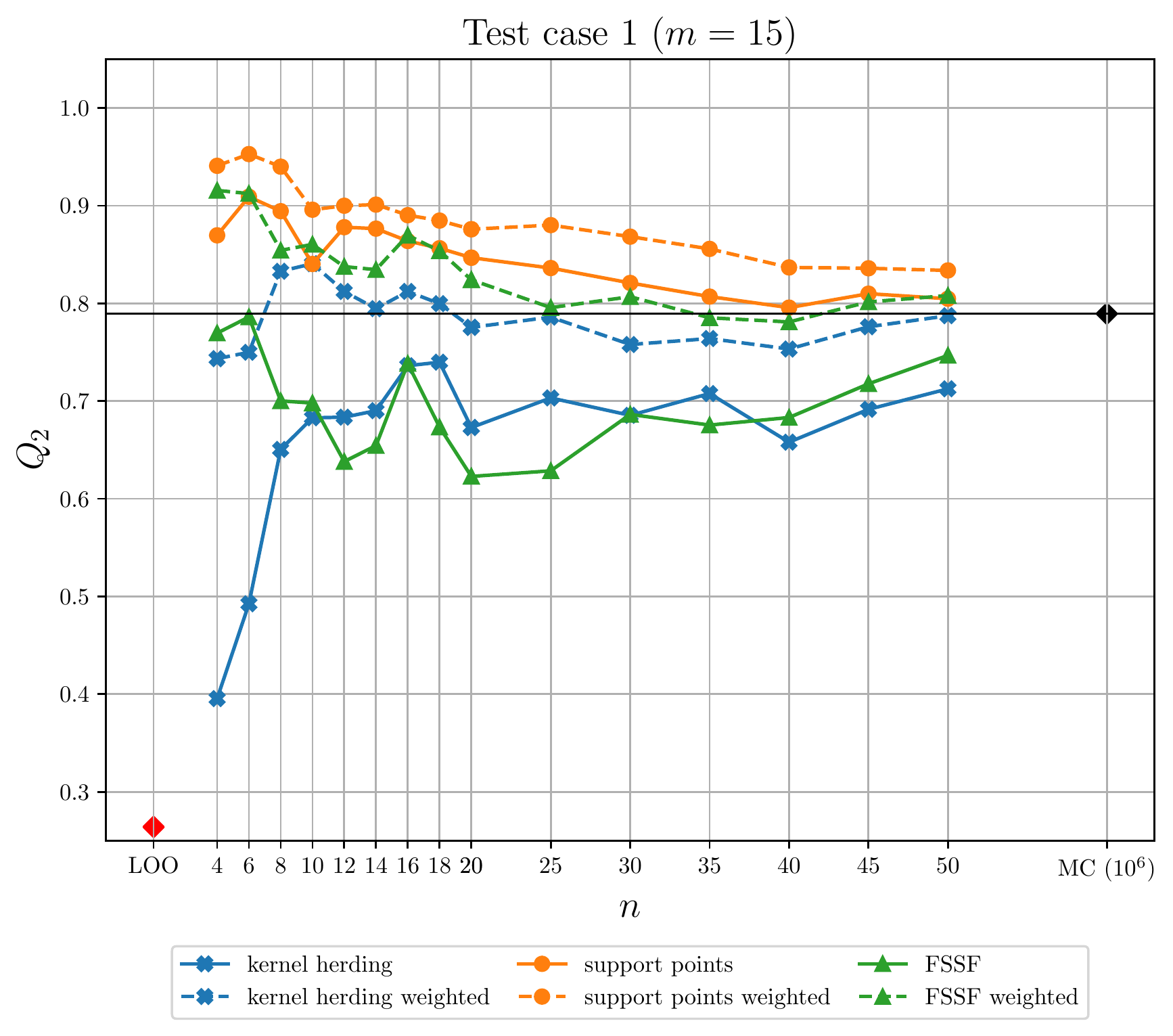}
        \quad
        \includegraphics[width=0.3\textwidth]{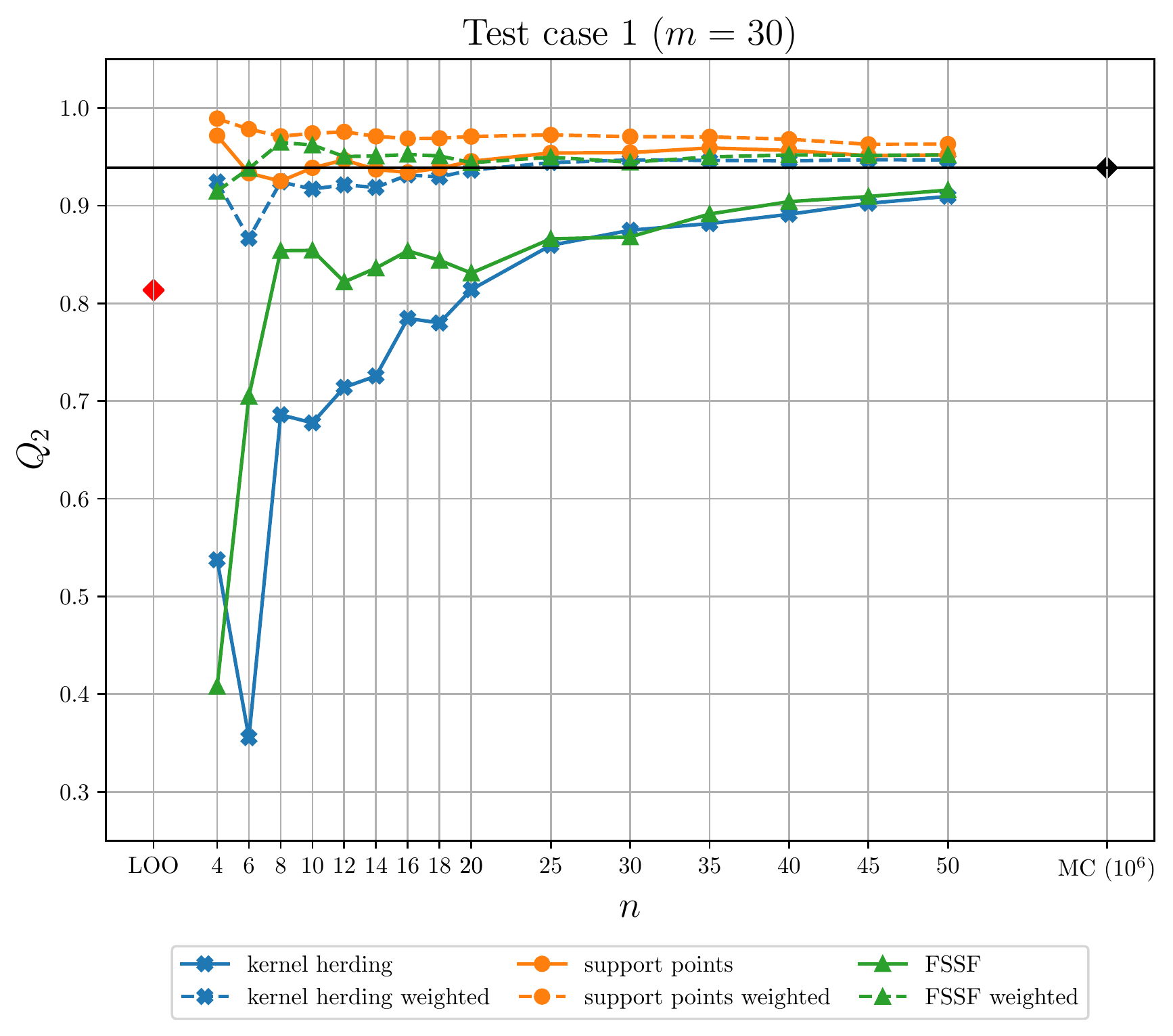}
        \caption{Test-case 1: predictivity assessment of a poor (left), good (center) and very good (right) model with kernel herding, support points and FSSF test sets.}
        \label{fig:irregular_benchmark}
    \end{minipage}

    \vspace{1cm}

    \begin{minipage}[t]{\textwidth}
        \includegraphics[width=0.3\textwidth]{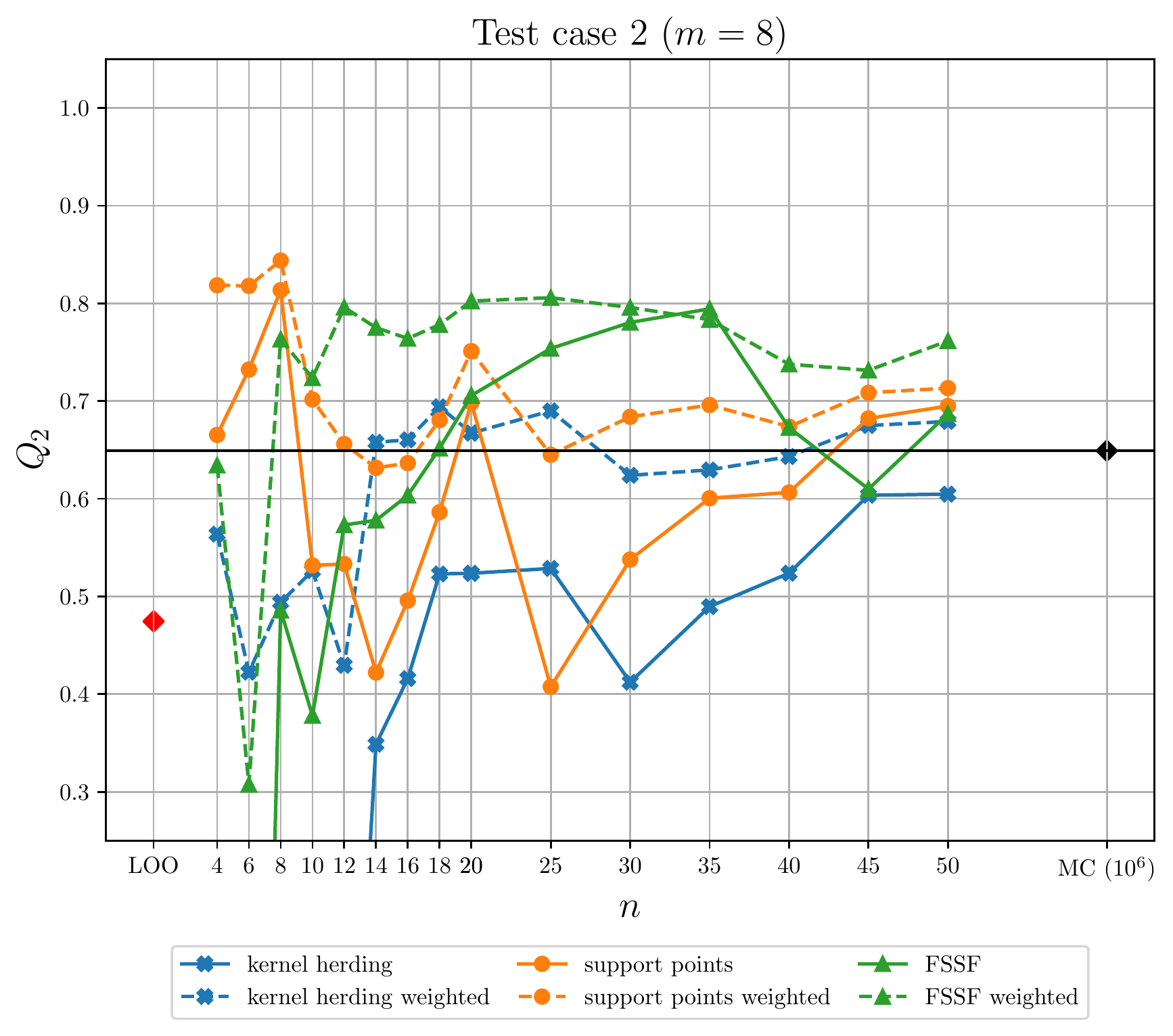}
        \quad
        \includegraphics[width=0.3\textwidth]{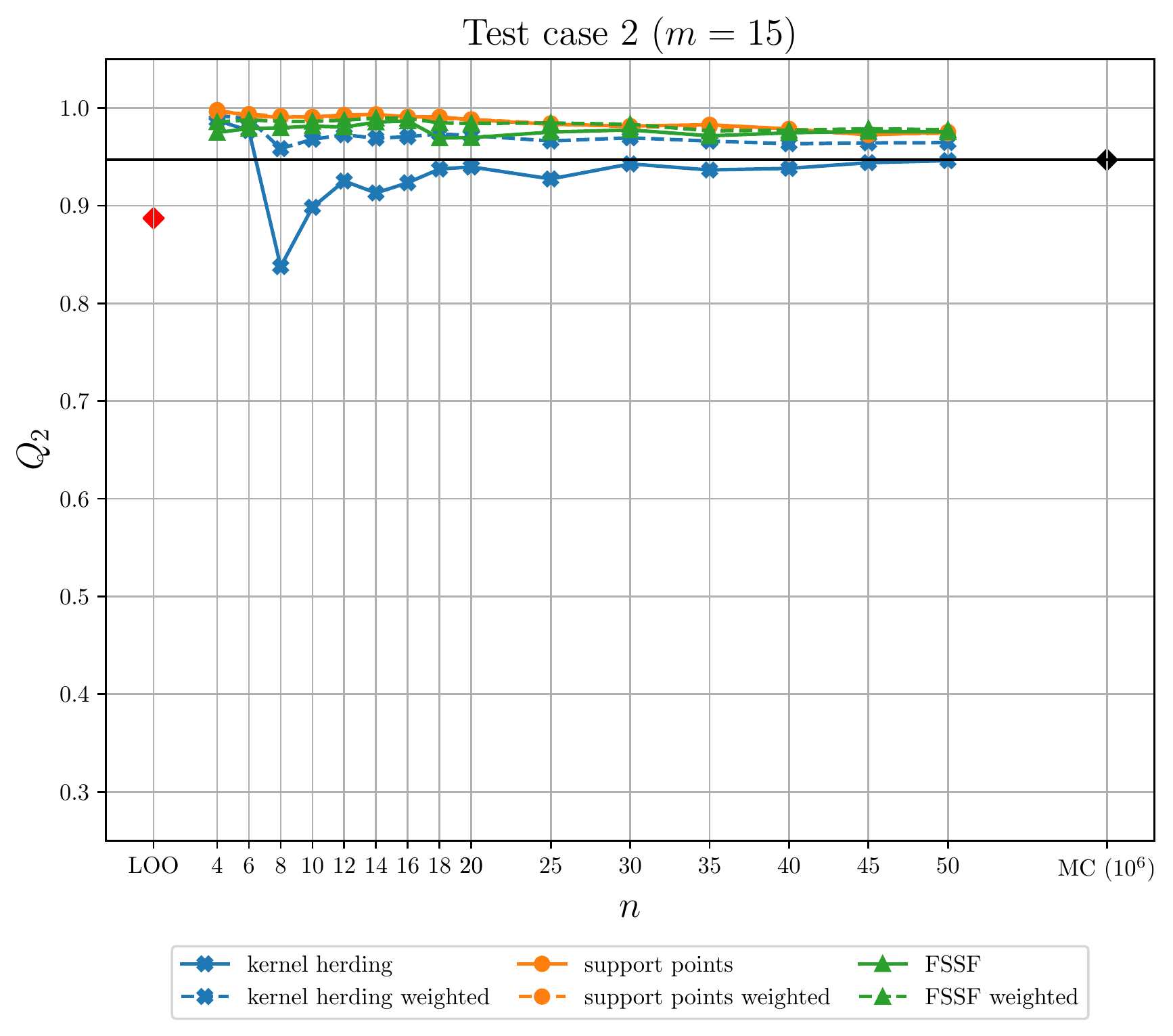}
        \quad
        \includegraphics[width=0.3\textwidth]{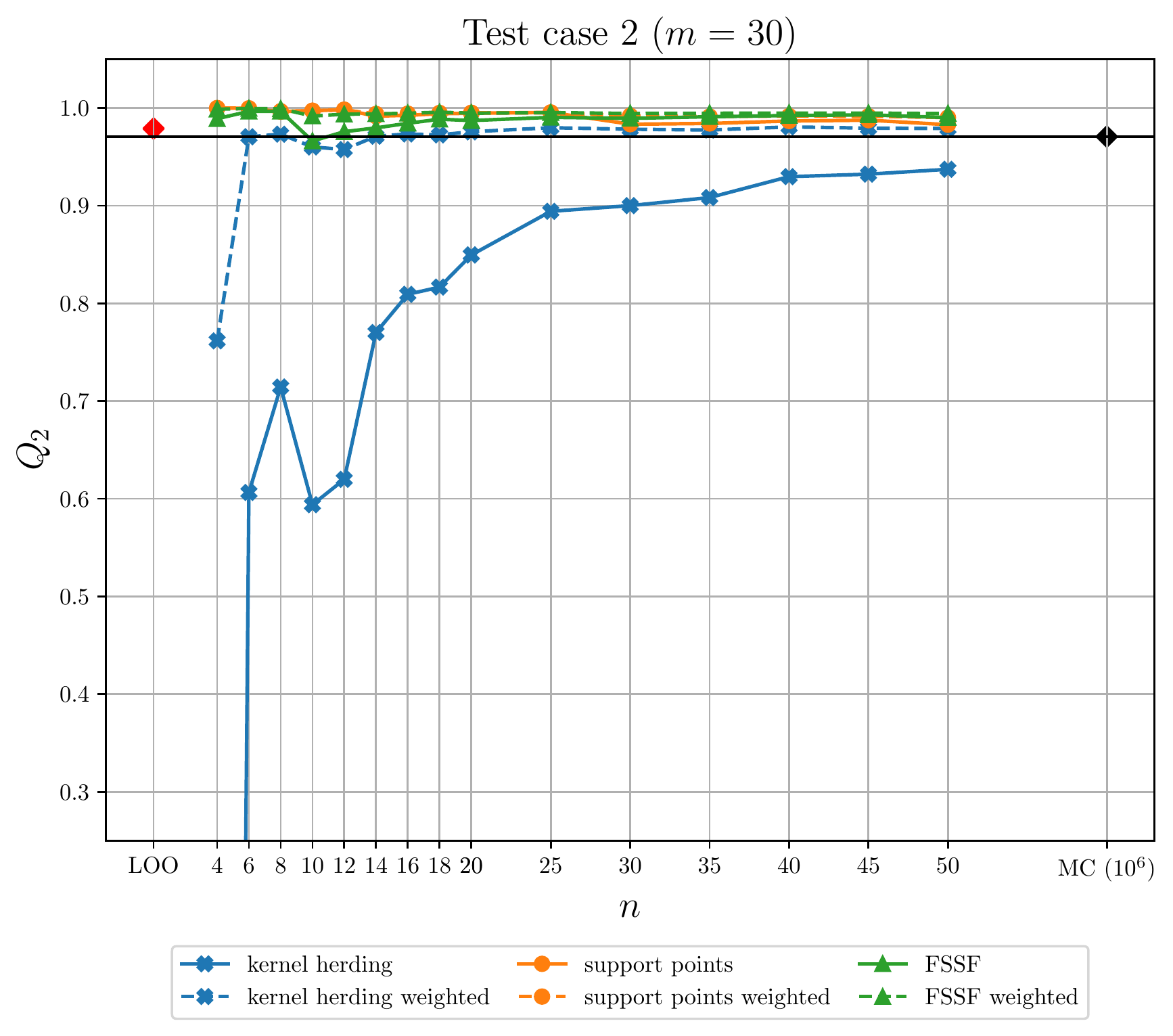}
        \caption{Test-case 2: predictivity assessment of a poor (left), good (center) and very good (right) model with kernel herding, support points and FSSF test sets.}
        \label{fig:cosin_benchmark}
    \end{minipage}

\end{sidewaysfigure}

\begin{sidewaysfigure}
    \vspace{11cm}
    \begin{minipage}[t]{\textwidth}
        \includegraphics[width=0.3\textwidth]{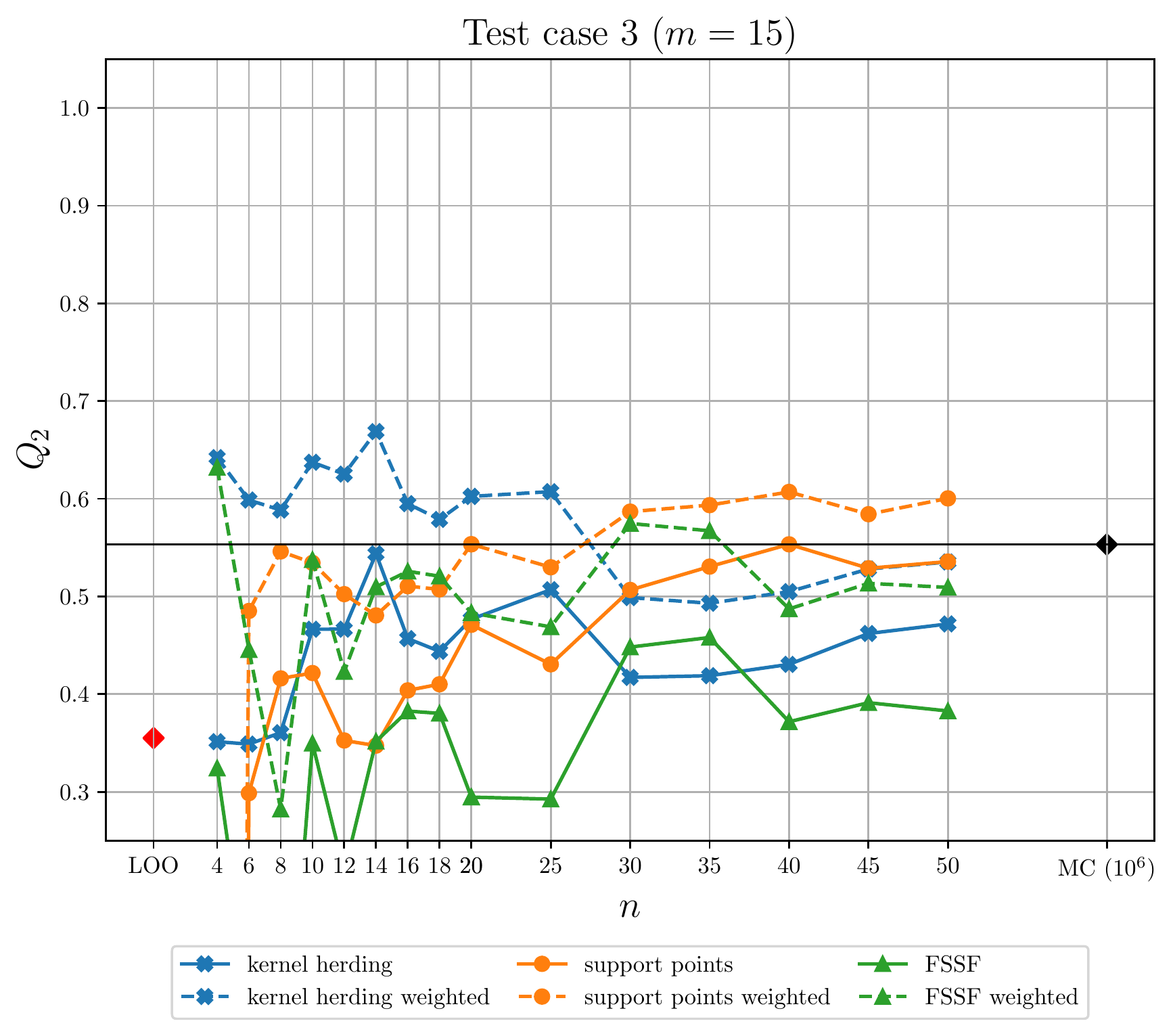}
        \quad
        \includegraphics[width=0.3\textwidth]{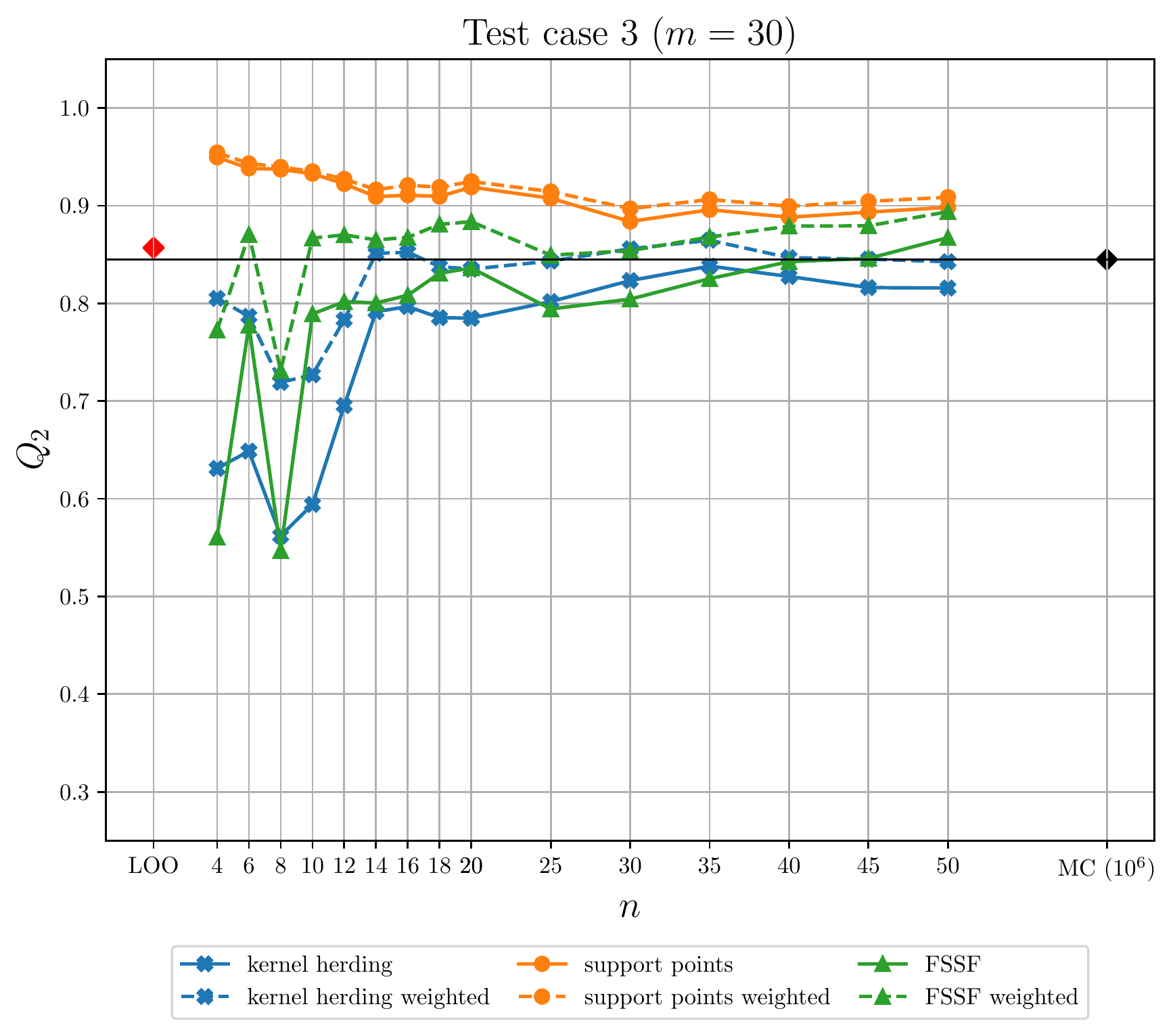}
        \quad
        \includegraphics[width=0.3\textwidth]{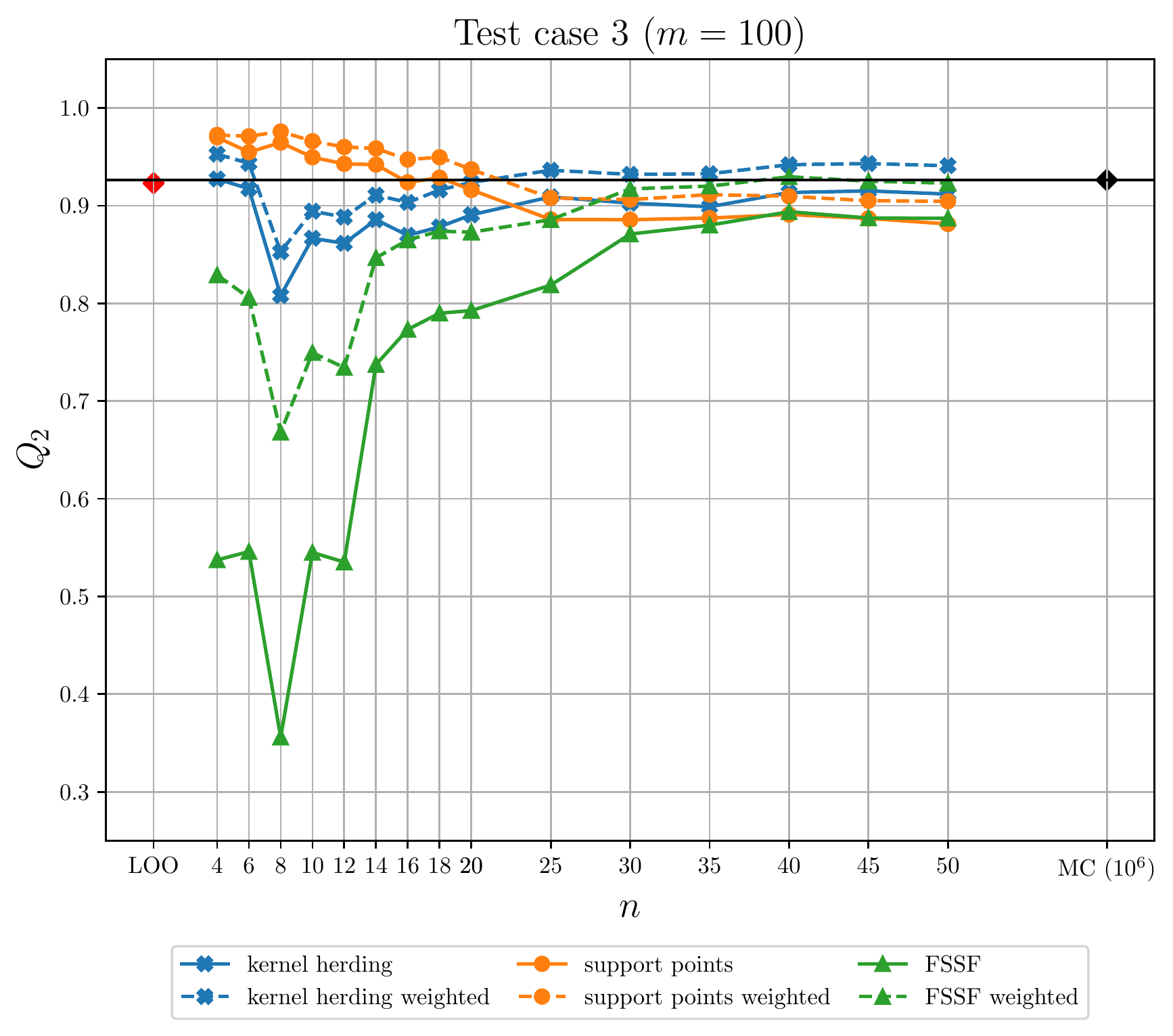}
        \caption{Test-case 3: predictivity assessment of a poor (left), good (center) and very good (right) model with kernel herding, support points and FSSF test sets.}
        \label{fig:gsobol_benchmark}
    \end{minipage}

\end{sidewaysfigure}

The numerical results obtained in this section are presented in Figures~\ref{fig:irregular_benchmark}, \ref{fig:cosin_benchmark}, and \ref{fig:gsobol_benchmark}. Each figure corresponds to one of the test-cases and gathers three sub-figures, corresponding to test sets with sizes $m$ yielding  poor (left), reasonably good (centre) or very good (right) fittings. 
 
The baseline value of $Q_{MC}^2$, calculated with $10^6$ Monte-Carlo points, is indicated by the black diamonds (the black  horizontal lines).
We assume that the error of $Q_{MC}^2$ is much smaller than the errors of all other estimators, and compare the distinct methods through their ability to approximate $Q_{MC}^2$. For each sequence of nested test-sets ($n\in\{4,\ldots,50\}$), the observed values of $\widehat Q^2_n$ (equation \eqref{eq:Q2test}) and $Q_{n*}^2$ (equation \eqref{Q2-est*}), are plotted as the solid and dashed lines, respectively.  

The figures also show the value $Q^2_{LOO}$ obtained by Leave-One-Out (LOO) cross validation, which is indicated at the left of each figure  by a red  diamond (values smaller than $0.25$ are not shown). 
Note that, contrarily to the other methods considered, for LOO the test set  is not disjoint from the training set, and thus the method does not satisfy the conditions set in the Introduction. As we repeat the complete model-fitting procedure for each training sample of size $m-1$, including the maximum-likelihood estimation of the correlation lengths of the Matérn kernel, the closed-form expressions of \cite{Dubrule83} cannot be used, making the computations rather intensive. 
As the three figures show, and as we should expect, $Q_{LOO}^2$ tends to under-estimate $Q_{\mathrm{ideal}}^2$:  by construction of the training set, LOO cross validation relies on  model predictions at points $\bx^{(i)}$ far from the other $m-1$ design points used to build the model, and thus tends to systematically overestimate the prediction error at $\bx^{(i)}$. The underestimation of $Q_{\mathrm{ideal}}^2$ can be particularly severe when $m$ is small, the training set being then necessarily sparse; see Figure~\ref{fig:irregular_benchmark} where $Q_{LOO}^2<0.3$ for $m=5$ and 15. 

Let us first concentrate on the non-weighted estimators (solid curves).
 We can see that the two MMD-based constructions,  support points (in orange) and kernel herding (in blue), generally produce better validation designs than FSSF (green curves), leading to values of $\widehat Q^2_n$ that approach $Q_{\mathrm{ideal}}^2$ quicker as $n$ increases. This is particularly noticeable for ``good'' and ``very good'' models (central and rightmost panels of all three figures).
This supports the idea that test sets should complement the training set $\bX_m$ by populating the holes it leaves in $\X$ while at the same time be able to mimic the target distribution $\mu$, this second objective being more difficult to achieve for FSSF than for the MMD-based constructions. 


Comparison of the two MMD based estimators reveals that support points tend to under-estimate ISE, leading to an over-confident assessment of the model predictivity, while kernel herding displays the expected behaviour, with a negative bias that decreases with $n$. The reason for the positive bias of estimates based on support points designs is not fully understood, but may be linked to the fact that support points tend to place validation points at ``mid-range'' from the designs (and not at the furthest points like FSSF or kernel herding), see central and rightmost panels in Figure \ref{fig:uniform_validation_designs}, and thus residuals at these points are themselves already better representatives of the local average errors. 

We consider now the impact of the GP-based weighting of the residuals when estimating $Q^2$ (by $Q_{n*}^2$), which is related to the relative training-set/validation-set geometry (the manner in which the two designs are entangled in  ambient space). The improvement resulting of applying residual weighting is apparent on all panels of the three figures, the dashed curves lying  closer to $Q_{\mathrm{ideal}}^2$ than their solid  counterparts;  see in particular kernel herding (blue curve) in Figure~\ref{fig:irregular_benchmark} and FSSF (green curve) in Figure~\ref{fig:cosin_benchmark}. Unexpectedly, the estimators based on support points seem to be rather insensitive to residual weighting, the dashed and solid orange curves being most of the time close to each other (and in any case, much closer that the green and blue ones). While the reason for this behavior deserves a deeper study, the fact that the support point designs -- see Figure \ref{fig:uniform_validation_designs} -- sample in a better manner the range of possible training-to-validation distances, being in some sense less space-filling than both FSSF and kernel herding, is again a plausible explanation for this weaker sensitivity to residual weighting.

Consider now comparison of the behaviour across test-cases. Setting aside the strikingly singular situation of test-case 2, for which kernel herding displays a pathological (bad) behaviour for the ``very good'' model, and all methods present an overall astonishing good  behaviour, we can conclude that the details of the tested function do not seem to play an important role concerning the relative merits of the estimators and validation designs.

We finally observe how the methods behave for models of distinct quality ($m$ leading to poor, good or very good models), comparing the three panels in each figure.
On the left panels, $m$ is too small for the model $\eta_m$ to be accurate, and all methods and test-set sizes are able to detect this.
For models of practical interest (good and very good), the test sets generated with support points and kernel herding allow a reasonably accurate estimation of $Q^2$ with a few points. Note, incidentally, that except for test-case 2 (where the interplay with a non-uniform measure $\mu$ complicates the analysis), it is in general easier to estimate the quality of the very good model (right-most panel) than that of the good model (central panel), indicating that the expected complexity (the entropy) of the residual process should be a key factor determining how large the validation set must be. In particular, it may be that larger values of $m$ allow for smaller values of $n$.

\section{Numerical results II: splitting a dataset into a training set and a test set}
\label{sec:5}



In this section, we illustrate the performance of the different designs and estimators considered in this paper when applied in the context of an industrial application, to split a given dataset of size $N$ into training  and test sets, with $m$ and $n$ points respectively, $m+n=N$. In contrast with \cite{josvak22}, the observations $y(\bx^{(i)})$, $i=1,\ldots,N$, are not used in the splitting mechanism, meaning that it can be performed before the observations are collected and that there cannot be any selection bias related to observations
(indeed, the use of observation values in a MMD-based splitting criterion may favour the allocation of the most different observations to different sets, training versus validation).

A ML model is fitted to the training data, and the data collected on the test-set are used to assess the predictivity of the model. The influence of the ratio $r_n=n/N=1-m/N$ on the quality assessment is investigated. We also consider Random Cross-Validation (RCV), where $n$ points are chosen at random among the $N$ points of the dataset: for each $n$, there are $N \choose n$ possible choices, and we randomly select $R=1\,000$ designs among them.
We fit a model to each of the $m$-point complementary designs ($m=N-n$), which yields an empirical distribution of $Q^2$ values for each ratio $n/N$ considered. 

\subsection{Industrial test-case CATHARE}

\begin{figure}
    \centering
    \includegraphics[width=0.35\textwidth]{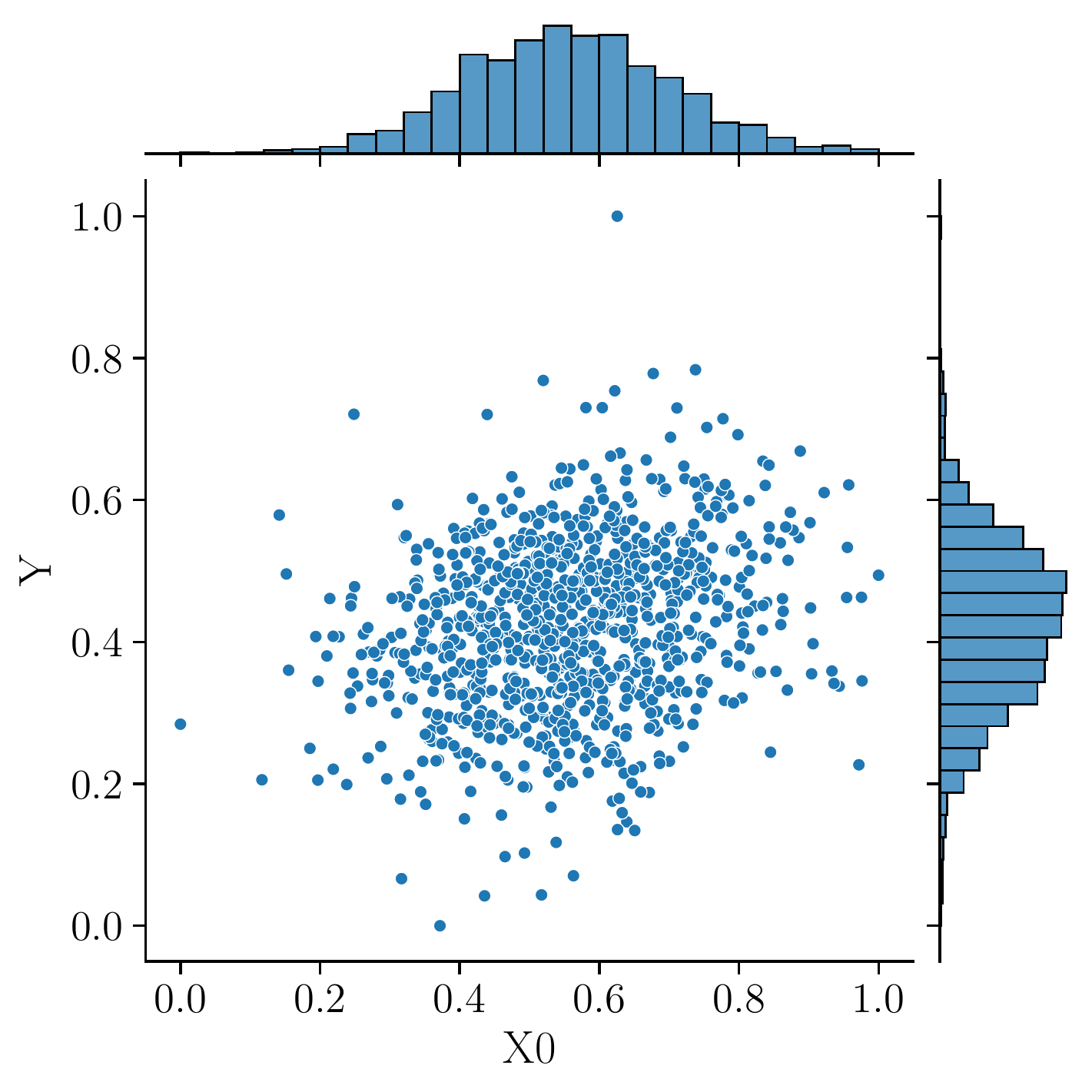}
    \includegraphics[width=0.35\textwidth]{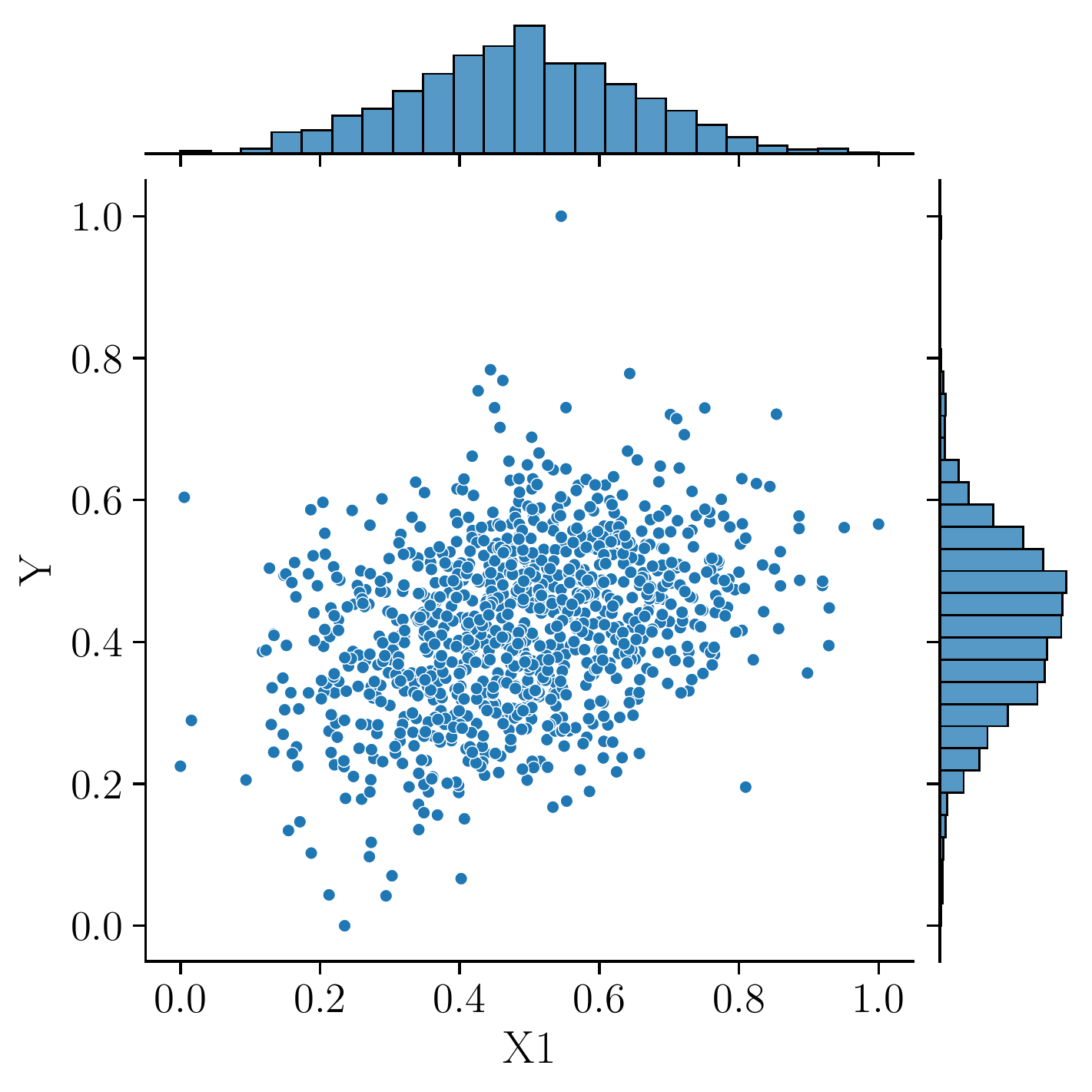}\\
    \includegraphics[width=0.35\textwidth]{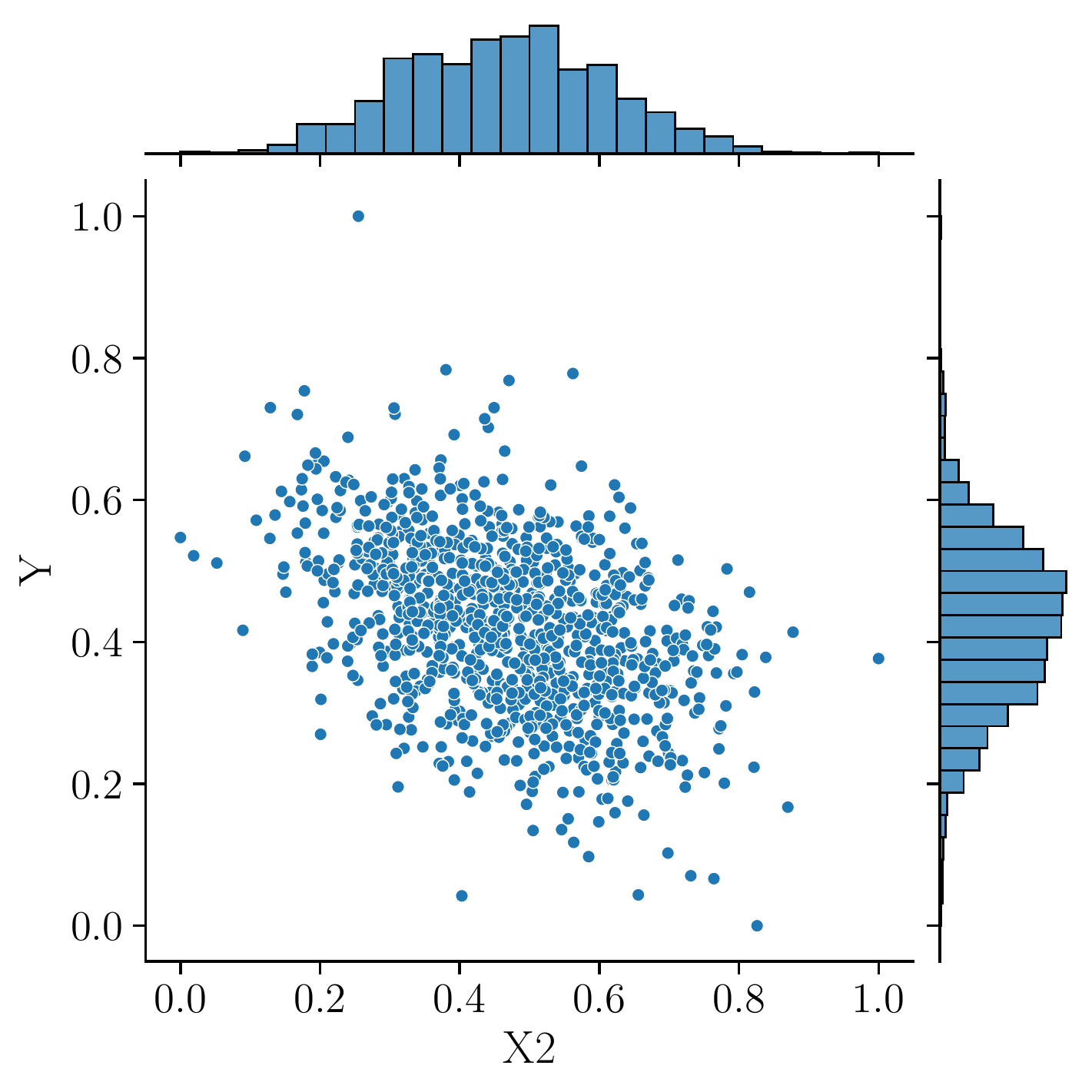}
    \includegraphics[width=0.35\textwidth]{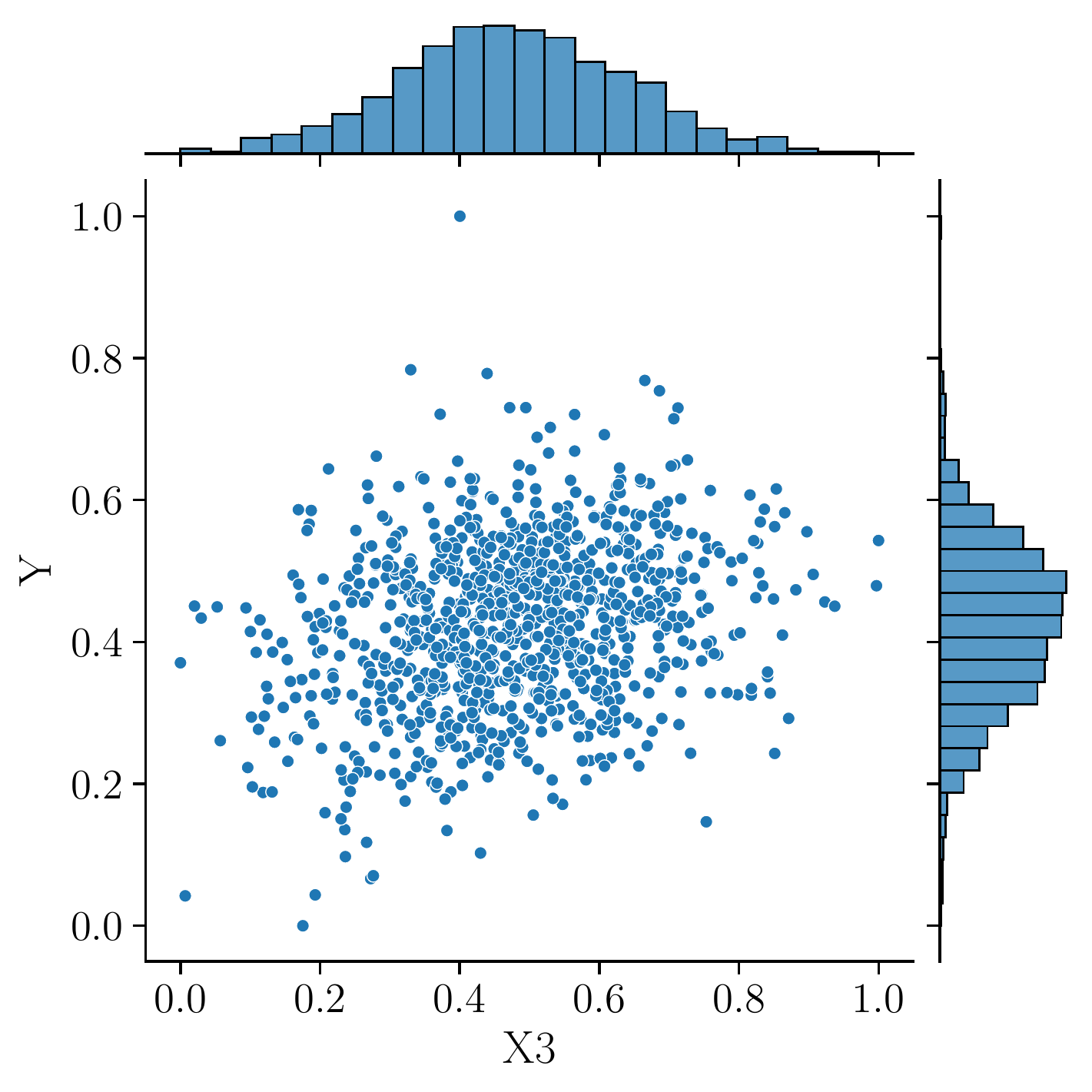}\\
    \includegraphics[width=0.35\textwidth]{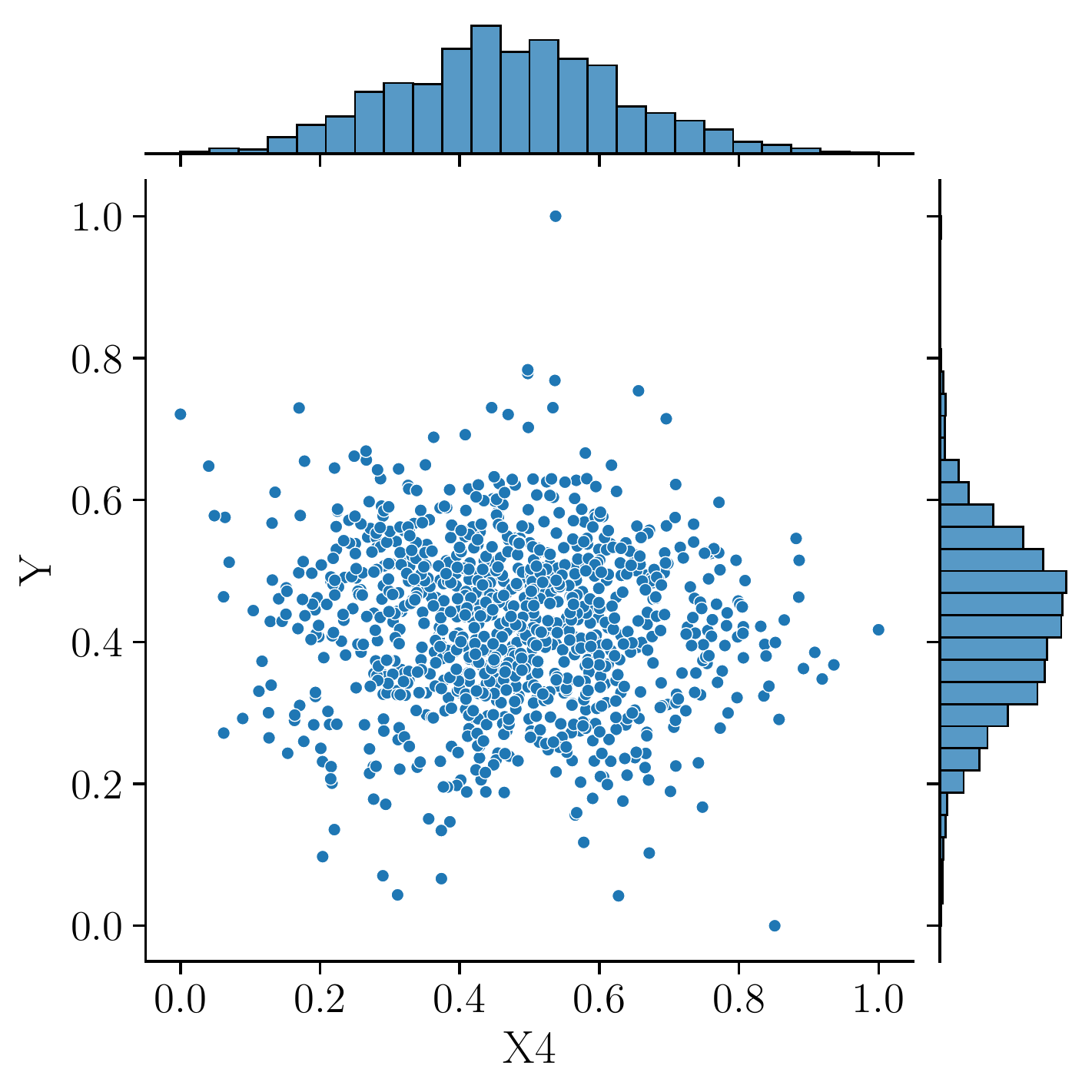}
    \includegraphics[width=0.35\textwidth]{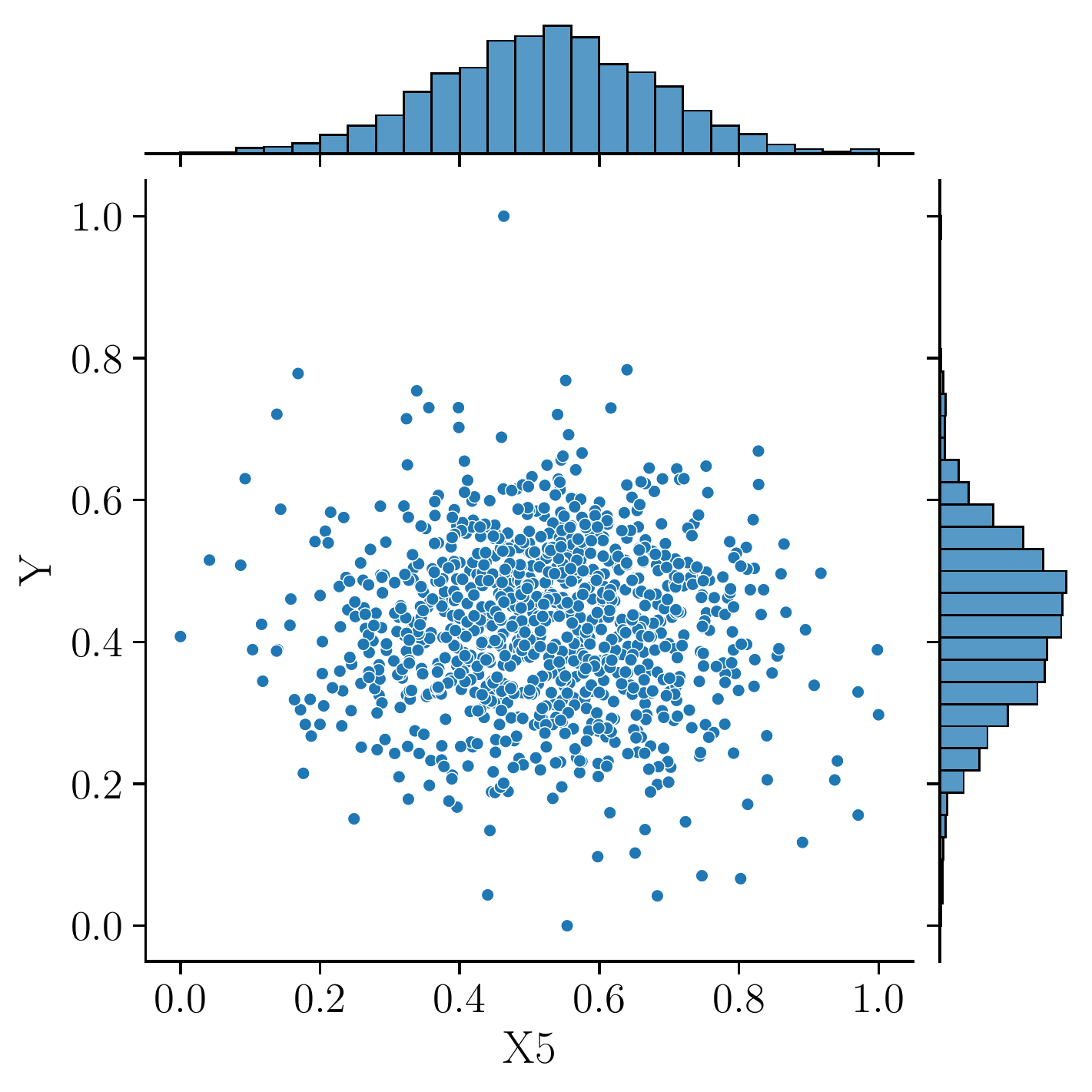}\\
    \includegraphics[width=0.35\textwidth]{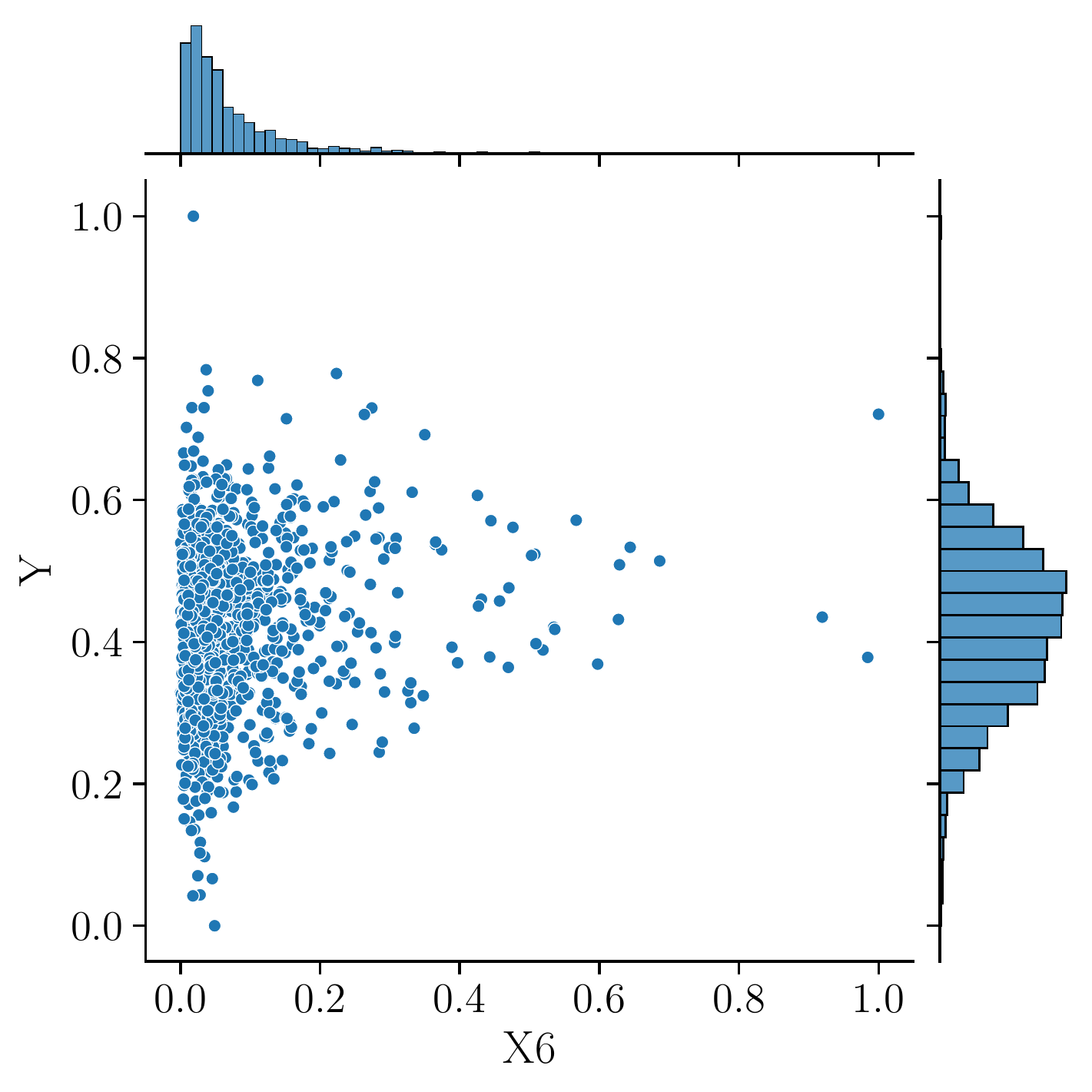}
    \includegraphics[width=0.35\textwidth]{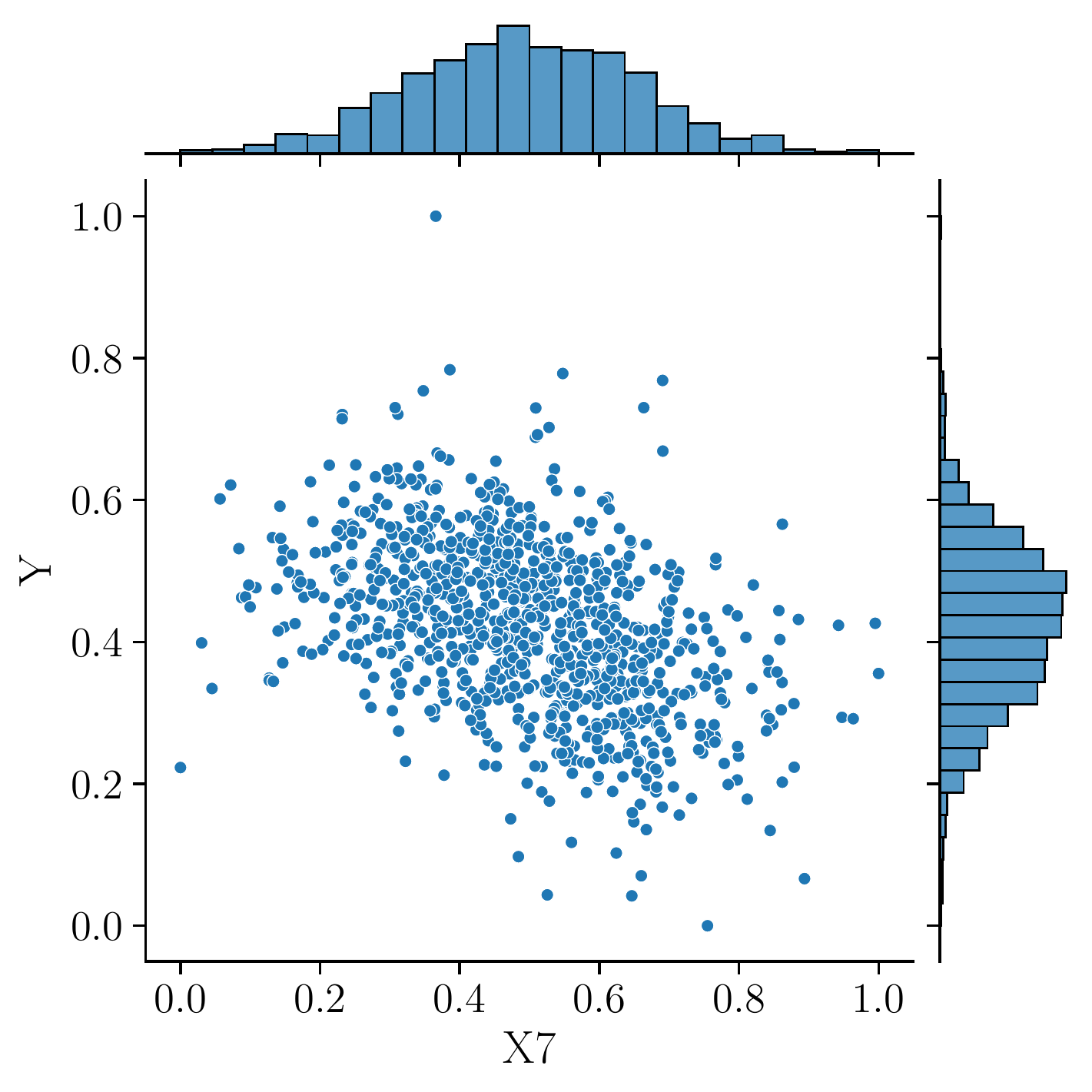}\\
    \includegraphics[width=0.35\textwidth]{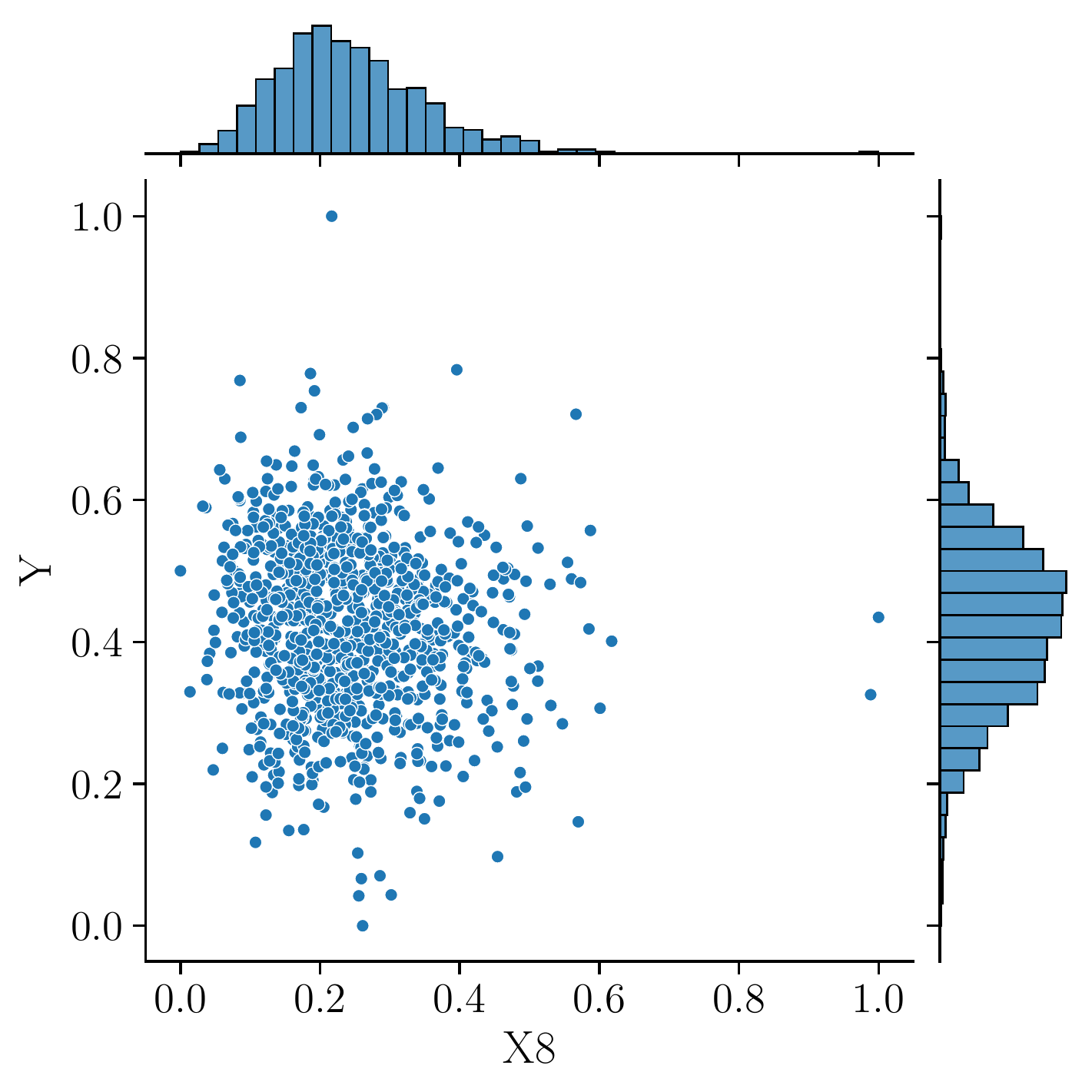}
    \includegraphics[width=0.35\textwidth]{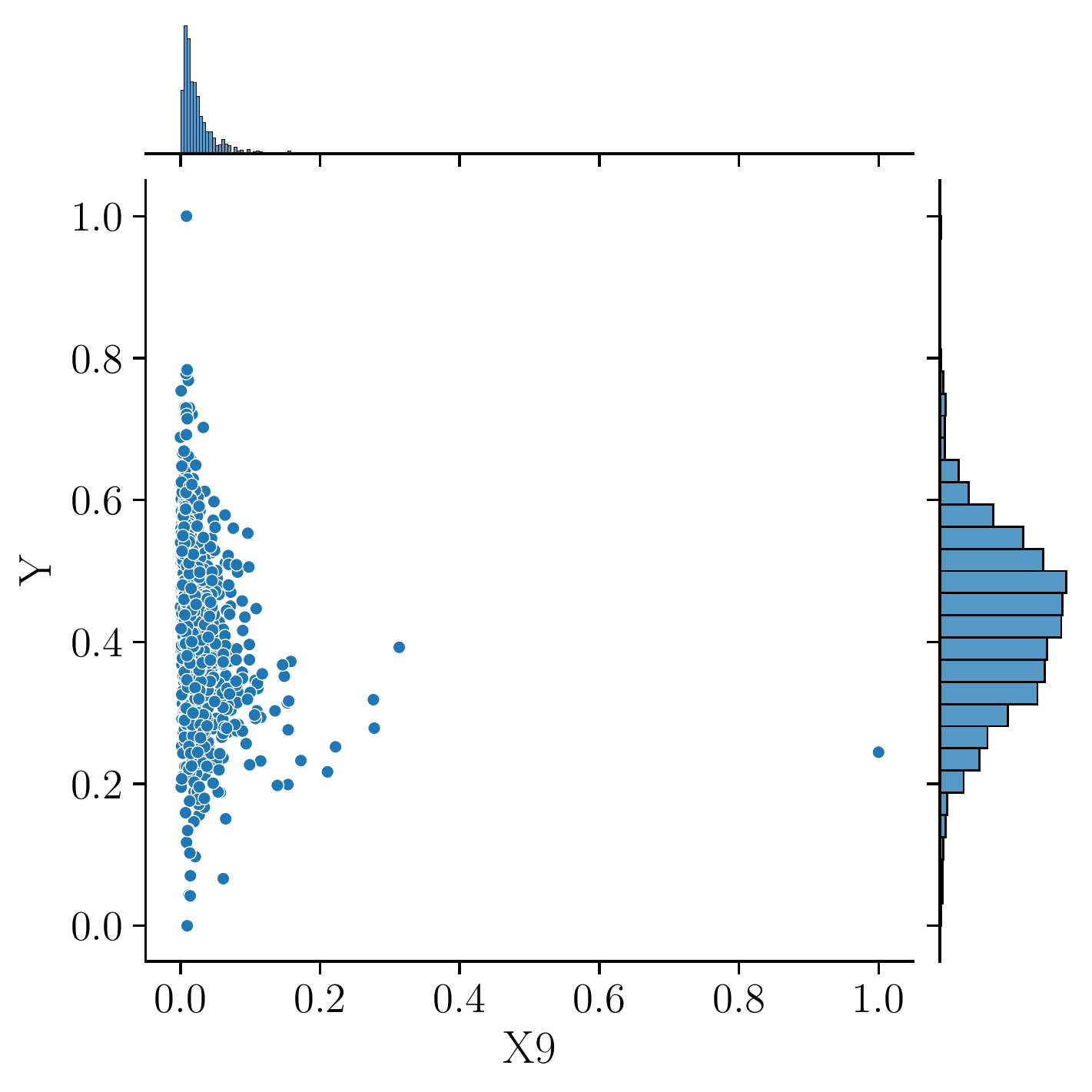}
    \caption{Test-case CATHARE: inputs output scatter plots ($N=10^3$) }
    \label{fig:cathare_paiplot}
\end{figure}

The test-case corresponds to  the computer code CATHARE2 (for ``Code Avanc\'e de ThermoHydraulique pour les Accidents de R\'eacteurs \`a Eau''), which models the thermal-hydraulic behavior inside nuclear pressurized water reactors \cite{gefant11}.
The studied scenario simulates a hypothetical large-break loss of primary coolant accident for which the output of interest is the peak cladding temperature \cite{decbaz08,ioobou10}.
The complexity of this application lies in the large run-time  of the computer model (of the order of twenty minutes) and in the high dimension of the input space: the model involves $53$  input parameters $z_i$, corresponding mostly to constants of physical laws, but also coding initial conditions, material properties and geometrical modeling. The $z_i$ were independently sampled according to normal or log-normal distributions (see axes histograms in Figure \ref{fig:cathare_paiplot} corresponding to $10$ inputs).
These characteristics make this test-case challenging in terms of construction of a surrogate model and validation of its predictivity.

%


    

We have access to an existing Monte Carlo sample $\bZ_N$ of $N=1\,000$ points in $\R^{53}$, that corresponds to $53$ independent random input configurations; see \cite{ioobou10} for details. The output  of  the  CATHARE2 code  at these $N$ points is also available. To reduce the dimensionality of this dataset, we first performed a sensitivity analysis \cite{davgam21} to eliminate inputs that do not impact the output significantly. This dimension-reduction step relies on the Hilbert-Schmidt Independence Criterion (HSIC), which is known as a powerful tool to perform input screening from a single sample of inputs and output values without reference to any specific ML regression model \cite{grebou05,dav15}.
HSIC-based statistical tests and their associated $p$-values are used to identify (with a $5\%$-threshold) inputs on which the output is significantly dependent (and therefore, also those of little influence). They were successfully applied to similar datasets from thermal-hydraulic applications in \cite{marcha21,marioo21}.
The screened dataset only includes $10$ influential inputs, over which the candidate set $\bX_N$ used for the construction of the test-set $\bX_n$ (and therefore of the complementary training set $\bX_{N-n}$) is defined. 
An input-output scatter plot is presented in Figure~\ref{fig:cathare_paiplot}, 
showing that indeed the retained factors are correlated with the code output. The marginal distributions are shown as histograms along to the axes of the plots.

To include RCV in the methods to be compared, we need to be able to construct many (here, $R=1\,000$) different models $\eta_m$ for each considered design size $m$.
Since Gaussian Process regression proved to be too expensive for this purpose, we settled for the  comparatively cheaper Partial Least Squares (PLS) method \cite{wolsjo01}, which retains acceptable accuracy.
For each given training set, the model obtained is a sum of monomials in the 10 input variables. Note that models constructed with different training sets may involve different monomials and have different numbers of monomial terms. 

\subsection{Benchmark results and analysis}

\begin{figure}
    \centering
    \includegraphics[width=0.9\textwidth]{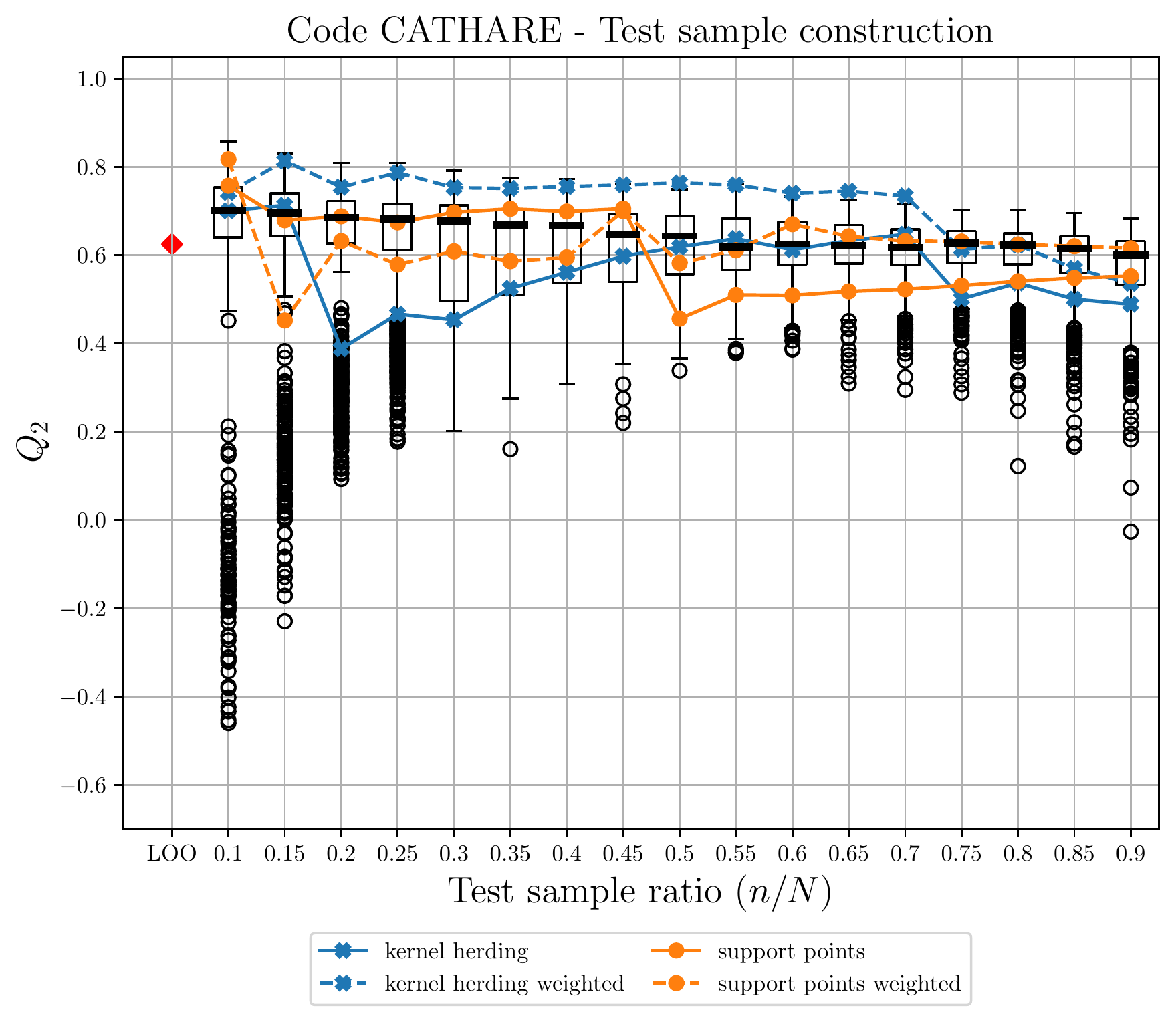}
    \caption{Test-case CATHARE: estimated $Q^2$. The box-plots are for random cross-validation, the red diamond (left) is for $Q^2_{LOO}$.}
    \label{fig:cathareC2_benchmark}
\end{figure}

Figure~\ref{fig:cathareC2_benchmark} compares various ways of extracting an $n$-point test set from an $N$-point dataset to estimate model predictivity, for different splitting ratios $n/N\in\{0.1,0.15,0.2,\ldots,0.9\}$. 

Consider RCV first. For each value of $r_n=n/N$, the empirical distribution of $Q^2_{RCV}$ obtained from $R=10^3$ random splittings of $\bX_N$ into $\bX_m\cup\bX_n$ is summarized by a boxplot. Depending on $r_n$, we can roughly distinguish three behaviors. For $0.1 \leq r_n \lesssim 0.3$ the distribution is bi-modal, with the lower mode corresponding to unlucky test-set selections leading to poor performance evaluations. When $0.3 \lesssim n/N \lesssim 0.7$, the distribution looks uni-modal, revealing a more stable performance evaluation. Note that this is (partly) in line with the recommendations discussed in section \ref{sec:1}. For $r_n \gtrsim  0.7$, the variance of the distribution increases with $r_n$: many unlucky  training sets lead to poor models. Note that the median of the empirical distribution slowly decreases as $r_n$ increases, which is consistent with the intuition that the model predictivity should decrease when the size of the training set decreases. 

For completeness, we also show  by a red diamond on the left of Figure~\ref{fig:cathareC2_benchmark} the value of $Q^2_{LOO}$ computed by LOO cross-validation. In principle, being computed using the entire dataset, this value should establish an upper bound on the quality of models computed with smaller training sets. 
This is indeed the case for small training sets (rightmost values in the figure), for which the predictivity estimated by LOO is above the majority of the predictivity indexes calculated. But at the same time, we know that LOO cross-validation tends to overestimate the errors, which explains the higher predictivity estimated by some other methods when $m=N-n$ is large enough.


Compare now the behavior of the two MMD-based algorithms of Section \ref{sec:2}, $\widehat Q^2_n$ (un-weighted) and $Q_{n*}^2$ (weighted) are plotted using solid and dashed lines, respectively, for both kernel herding (in blue) and support points (in orange). FSSF test-sets are not considered, as the application of an iso-probabilistic transformation imposes knowledge of the input distribution, which is not known for this example. 
Compare first the unweighted versions of the two MMD-based estimators. 
For small values of the ratio $r_n$, $0.1 \lesssim r_n \lesssim 0.45$, the relative behavior of support points and kernel herding coincides with what we observed in the previous section, support points (solid orange line) estimating a better performance than kernel herding (solid blue line), which, moreover, is close to the median of the empirical distribution of $Q^2_{RCV}$. However, for $r_n \geq 0.5$, the dominance is reversed, support points estimating a worse performance than kernel herding. 

As $r_n$ increases up to $r_n \lesssim 0.7$ the solid orange and blue curves crossover, and it is now $\widehat Q^2_n$ for kernel herding that approximates the RCV empirical median, while the value obtained with support points underestimates the predictivity index. Also, note that for (irrealistic) very large values of $r_n$ both support points and kernel herding estimate lower $Q^2$ values, which are smaller than the median of the RCV estimates.

Let us now focus on the effect of residual weighting, i.e., in estimators $Q_{n*}^2$ which use the weights computed by the method of Section~\ref{S:weighting}, shown in dashed lines in Figure \ref{fig:cathareC2_benchmark}. First, note that while for kernel herding weighting leads, as in the previous section, to higher estimates of the predictivity (compare solid and dashed blue lines), this is not the case for support points (solid and dashed orange curves), which, for small split ratios, produces smaller estimates when weighting is introduced. In the large $r_n$ region, the behavior is consistent with what we saw previously, weighting inducing an increase of the estimated predictivity. It is remarkable -- and rather surprising -- that $Q_{n*}^2$ for support points (the dashed orange line) does not present the discontinuity of the uncorrected curve. 

\begin{figure}
    \centering
    \includegraphics[width=0.8\textwidth]{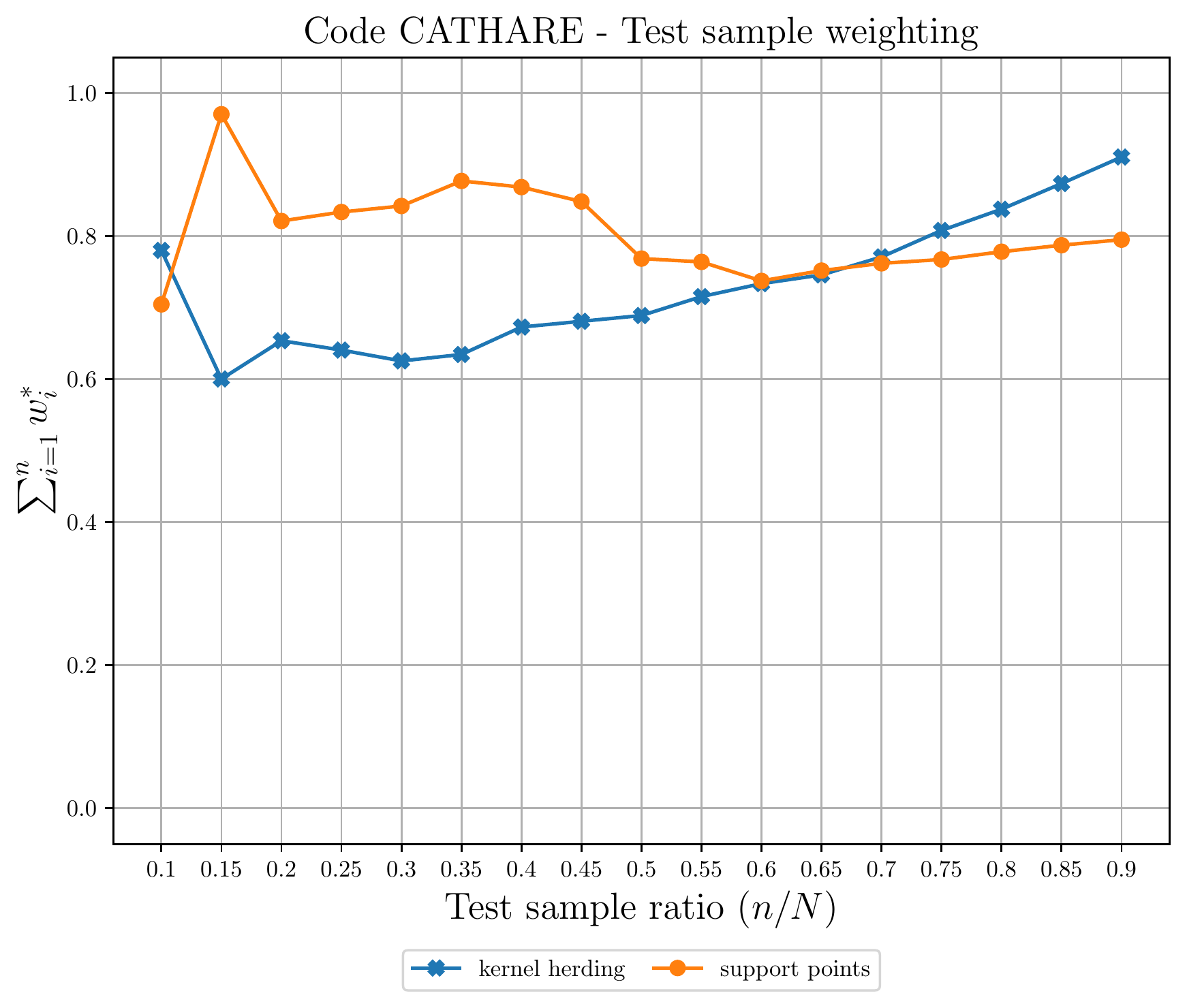}
    \caption{Test-case CATHARE: sum of the weights \eqref{Eq:weights}.}
    \label{fig:catharec2_weights}
\end{figure}

The sum $\sum_{i=1}^n w_i^*$ of the optimal weights of support points and kernel herding \eqref{Eq:weights} is shown in  Figure~\ref{fig:catharec2_weights} (orange and blue curves, respectively). The slow increase with $n/N$ of the sum of kernel-herding weights (blue line) is consistent with the increase of the volume of the input region around each validation point when the size of the training set decreases. 
The behavior of the sum of weights is more difficult to interpret for support points (orange line) but is consistent with the behavior of $Q_{n*}^2$ on Figure~\ref{fig:cathareC2_benchmark}. Note that the energy-distance kernel \eqref{eq:kE} used for support points cannot be used for the weighting method of Section~\ref{S:weighting} as $K_E$ is not positive definite but only conditionally positive definite. A full understanding of the observed curves would require a deeper analysis of the geometric characteristics of the designs generated by the two MMD methods, in particular of their interleaving with the training designs, which is not compatible with the space constraints of this manuscript.

While a number of unanswered points remain, in particular how deeply the behaviours observed may be affected by the poor predictivity resulting from the chosen PLS modeling methodology, the example presented in this section shows that  the construction of test sets via MMD minimization and estimation of the predictivity index using the weighted estimator $Q_{n*}^2$ is promising as an efficient alternative to RCV: at a much lower computational cost, it  builds performance estimates based on independent data the model developers may not have access to. Moreover, kernel herding proved, in the examples studied in this manuscript, to be a more reliable option for designing the test set, exhibiting a behavior that is consistent with what is expected, and very good estimation quality when the residuals over the design points are appropriately weighted.



\vspace{-1cm}

\section{Conclusion}
\label{sec:6}


Our study shows that ideas and tools from the design of experiments framework can be transposed to the problem of test-set selection. This paper explored approaches based on support points, kernel herding and FSSF, considering the incremental construction of a test set (\textit{i}) either as a particular space-filling design problem, where design points should populate the holes left in the design space by the training set, or (\textit{i})  from the point of view of partitioning a given dataset into a training set and a test set. 

A numerical benchmark has been performed for a panel of test-cases of different dimensions and complexity. Additionally to the usual predictivity coefficient, a new weighted metric (see \cite{PR2021a}) has been proposed and shown to improve assessment of the predictivity  of a given model for a given test set. 

This weighting procedure appears very efficient for interpolators, like Gaussian process regression models, as it corrects the bias when the points in the test set used to predict the errors are far from the training points. For the first three test-cases (Section~\ref{sec:4}), pairing one iterative design method with the weight-corrected estimator of the predictivity coefficient $Q^2$ shows promising results as the estimated $Q^2$ characteristic is close to the true one even for  test-sets of moderate size. 

Weighting can also be applied to models that do not interpolate the training data. For the industrial test-case of Section~\ref{sec:5}, the true $Q^2$ value is unknown, but the weight-corrected estimation $Q_{n*}^2$ of $Q^2$ is close to the value estimated by Leave-One-Out cross validation and to the median of the empirical distribution of $Q^2$ values obtained by random $k$-fold cross-validation. At the same time, estimation by $Q_{n*}^2$ involves a much smaller computational cost than cross-validation methods, and uses a dataset fully independent from the one used to construct the model.  

To each of the design methods considered to select a test set a downside can be attached. FSSF requires knowledge of the input distribution to be able to apply an iso-probabilistic transformation if necessary; it tends to select many points along the boundary of the candidate set considered. 
Support points require the computation of the $N(N-1)/2$ distances between all pairs of candidate points, which implies important memory requirements for large $N$; the energy-distance kernel on which the method relies cannot be used for the weighting procedure. Finally, the efficient implementation of kernel herding relies on  analytical expressions for the potentials $P_{K,\mu}$, see Appendices A and B, which are available for particular distributions (like the uniform and the normal) and kernels (like Matérn) only. The great freedom in the choice of the kernel $K$ gives a lot of flexibility, but at the same time implies that some non-trivial decisions have to be made; also, the internal parameters of $K$, such as its correlation lengths, must to be specified. Future work should go beyond empirical rules of thumb and study the influence of these choices.

We have only computed numerical tests with independent inputs. Kernel herding and support points are both well suited for probability measures not being equal to the product of their marginals, which is a frequent case in real datasets.
We have also only considered incremental constructions, as they allow to stop the validation procedure as soon as the estimation of the model predictivity is deemed sufficiently accurate, but it is also possible to select several points at once, using support points \cite{makjos18}, or MMD minimization in general \cite{TeymurGRO2021}. 

Further developments around this work could be as follows. Firstly, the incremental construction of a test set could be coupled with the definition of an appropriate stopping rule, in order to decide when  it is  necessary to continue improving the model (possibly by supplementing the initial design with the test set, which seems well suited to this). The MMD $d_{\Kbarbar}(\zeta_n^*,\mu)$ of Section~\ref{S:weighting} could play an important role in the derivation of such a rule. Secondly, the approach presented gives equal importance to all the $d$ inputs. However, it seems that inputs with a negligible influence on the output should receive less attention when selecting a test set. A preliminary screening step that identifies the important inputs would allow the test-set selection algorithm to be applied on these variables only. 
For example, when a $\bX_N\subset\R^d$ dataset is to be partitioned into $\bX_m\cup\bX_n$, one could use only $d'<d$ components to define the partition, but still use all $d$ components to build the model and estimate its (weighted) $Q^2$. 
Note, however, that this would imply a slight violation of the conditions mentioned in introduction, as it renders the test set dependent on the function observations.

Finally, in some cases the probability measure $\mu$ is known up to a normalizing constant. The use of a Stein kernel then makes the potential $P_{K,\mu}$ identically zero \cite{ChenMGBO2018,ChenBBGGMO2019}, which would facilitate the application of kernel herding. Also, more complex problems involve functional inputs, like temporal signals or images, or categorical variables; the application of the methods presented to kernels specifically designed for such situations raises challenging issues.


\begin{acknowledgement}
This work was supported by project INDEX (INcremental Design of EXperiments) ANR-18-CE91-0007 of the French National Research Agency (ANR).
The authors are grateful to Guillaume Levillain and Thomas Bittar for their code development during their work at EDF.
Thanks also to S\'ebastien Da Veiga for fruitful discussions.

\end{acknowledgement}


\section*{Appendix}
\addcontentsline{toc}{section}{Appendix}

\subsection*{Appendix A: Maximum Mean Discrepancy}\label{sec:MMD}

Let $K$ be a positive definite kernel on $\X\times\X$, defining a reproducing kernel Hilbert space (RKHS) $\HH_K$ of functions on $\X$, with scalar product $\langle f,g\rangle_{\HH_K}$ and norm $\|f\|_{\HH_K}$; see, e.g., \cite{bertho04}. For any $f\in\HH_K$ and any probability measures $\mu$ and $\xi$ on $\X$, we have
\begin{eqnarray}
\left| \int_\X f(\bx)\, \dd\xi(\bx) - \int_\X f(\bx)\, \dd\mu(\bx) \right| &=& \left| \int_\X \langle f,K_\bx\rangle_{\HH_K}\, \dd(\xi-\mu)(\bx) \right| \nonumber \\
&=& \left| \langle f,(P_{K,\xi}-P_{K,\mu}\rangle_{\HH_K} \right| \,, \label{eq:diff-of-int}
\end{eqnarray}
where we have denoted $K_\bx(\cdot)=K(\bx,\cdot)$ and used the reproducing property $f(\bx)=\langle f,K_\bx\rangle_{\HH_K}$ for all $\bx\in\X$, and where, for any probability measure $\nu$ on $\X$ and $\bx\in\X$,
\begin{equation}\label{eq:potnu}
P_{K,\nu}(\bx) = \int_\X K(\bx, \bx') \, \dd\nu(\bx') \,,
\end{equation}
is the potential of $\nu$ at $\bx$. $P_{K,\nu}\in\HH_K$ and is called kernel embedding of $\nu$ in ML.
In some cases, the potential can be expressed analytically (see. Appendix \ref{sec:potential}), otherwise it can be estimated by numerical quadrature (Quasi Monte Carlo). 
Cauchy-Schwartz inequality applied to \eqref{eq:diff-of-int} gives
$$
\left| \int_\X f(\bx)\, \dd\xi(\bx) - \int_\X f(\bx)\, \dd\mu(\bx) \right| \leq \|f\|_{\HH_K}\,\|P_{K,\xi}-P_{K,\mu}\|_{\HH_K}
$$
and therefore 
$$
\|P_{K,\xi}-P_{K,\mu}\|_{\HH_K}=\sup_{f\in\HH_K:\ \|f\|_{\HH_K}=1} \left| \int_\X f(\bx)\, \dd\xi(\bx) - \int_\X f(\bx)\, \dd\mu(\bx) \right|\,.
$$
The Maximum Mean Discrepancy (MMD) between $\xi$ and $\mu$ (for the kernel $K$ and set $\X$) is $d_K(\xi,\mu)=\|P_{K,\xi}-P_{K,\mu}\|_{\HH_K}$. Direct calculation gives 
\begin{eqnarray}
d_K^2(\xi,\mu) &=& \|P_{K,\xi}-P_{K,\mu}\|_{\HH_K}^2 = \int_{\X^2} K(\bx,\bx')\, \dd(\xi-\mu)(\bx)\dd(\xi-\mu)(\bx') \label{eq:MMD1} \\
&=& \EE_{\zeta,\zeta'\sim \xi} K(\zeta,\zeta') + \EE_{\zeta,\zeta'\sim \mu} K(\zeta,\zeta') - 2\EE_{\zeta\sim \xi, \zeta'\sim \mu} K(\zeta,\zeta') \,, \label{eq:MMD}
\end{eqnarray}
where the random variables $\zeta$ and $\zeta'$ in \eqref{eq:MMD} are independent, see \cite{smogre07}. When $K$ is the energy distance kernel \eqref{eq:kE}, one recovers the expression \eqref{eq:energySP} for the corresponding MMD. 
One may refer to \cite{srigre10} for an illuminating exposition on MMD, kernel embedding, and conditions on $K$ (the notion of characteristic kernel) that make $d_K$ a metric on the space of probability measures on $\X$. The distance and Mat\'ern kernels considered in this paper are characteristic. 

\subsection*{Appendix B: Analytical computation of potentials for Mat\'ern kernels} \label{sec:potential}

As for tensor-product kernels, the potential is the product of the one-dimensional potentials, we only consider one-dimensional input spaces.

For $\mu$ the uniform distribution on $[0,1]$ and $K$ the Matérn kernel $K_{5/2,\theta}$ with smoothness $\nu=5/2$ and correlation length $\theta$, see \eqref{eq:Matern5/2}, we get
\bea
    P_{K_{5/2,\theta},\mu}(x) = \frac{16 \theta}{3 \sqrt{5}} - \frac{1}{15 \theta} (S_\theta(x) + S_\theta(1-x)),
\eea
where
\bea
\nonumber
    S_\theta(x) = 
    \exp\left(- \frac{\sqrt{5}}{\theta} x \right) 
    \left( 5 \sqrt{5} x^2 + 25 \theta x + 8 \sqrt{5} \theta^2 \right).
\eea
The expressions $P_{K_{\nu,\theta},\mu}(x)$ for $\nu=1/2$ and $\nu=3/2$ can be found in \cite{prozhi20}. 

When $\mu$ is the standard normal distribution $\mathcal{N}(0,1)$, the potential $P_{K_{5/2,\theta},\mathcal{N}(0,1)}$ is 
$    P_{K_{5/2,\theta},\mathcal{N}(0,1)}(x) = T_\theta(x) + T_\theta(-x),
$
where
\bea
T_\theta(x) &=&
    \frac{1}{6} 
    \left( 
        \frac{5}{\theta^2} x^2 + 
        \left( 3 - \frac{10}{\theta^2} \right) \frac{\sqrt{5}}{\theta} x + \frac{5}{\theta^2} \left( \frac{5}{\theta^2} -2 \right) + 3
    \right) \\
    && \times\, \mathrm{erfc} \left( \frac{\frac{\sqrt{5}}{\theta} - x}{\sqrt{2}} \right)
    \exp\left(\frac{5}{2 \theta^2} - \frac{\sqrt{5}}{\theta}x\right) + \frac{1}{3 \sqrt{2 \pi}} \frac{\sqrt{5}}{\theta} 
    \left( 3 - \frac{5}{\theta^2} \right) \exp\left(-\frac{x^2}{2}\right). 
\eea


\end{document}